\documentclass{amsart}

\usepackage{amssymb,latexsym,amsfonts,amsmath}
\usepackage{graphicx}
\usepackage{diagrams}
\usepackage{tikz}
\usetikzlibrary{automata}

\usepackage{amssymb}
\usepackage{amsfonts}
\usepackage{amssymb,latexsym,amsfonts,amsmath}
\usepackage{graphicx}

\newtheorem{theorem}{Theorem}[section]
\newtheorem{lemma}[theorem]{Lemma}

\newtheorem{proposition}[theorem]{Proposition}
\newtheorem{corollary}[theorem]{Corollary}

\theoremstyle{definition}
\newtheorem{definition}[theorem]{Definition}
\newtheorem{example}[theorem]{Example}

\theoremstyle{remark}

\numberwithin{equation}{section}

\begin{document}
\begin{abstract}
Symbolic models are abstract descriptions of continuous systems in which symbols represent aggregates of continuous states. 
In the last few years there has been a growing interest in the use of symbolic models as a tool for mitigating complexity in control design. In fact, symbolic models
%
enable the use of well known algorithms in the context of supervisory control and algorithmic game theory, for controller synthesis.
%
%
Since the 1990's many researchers faced the problem of identifying classes of dynamical and control systems that admit symbolic models.
In this paper we make a further progress along this research line by 
focusing on control systems affected by disturbances. Our main contribution is to show that incrementally globally asymptotically stable nonlinear control systems with disturbances admit symbolic models. When specializing these results to linear systems, we show that these symbolic models can be easily constructed.
\end{abstract}

\title[Symbolic Models for Nonlinear Control Systems: Alternating Approximate Bisimulations]{Symbolic Models for Nonlinear Control Systems:\\ Alternating Approximate Bisimulations}
\thanks{This work has been partially supported by the National Science Foundation CAREER award 0717188.}

\author[Giordano Pola and Paulo Tabuada]{Giordano Pola and Paulo Tabuada}
\address{Department of Electrical Engineering\\
University of California at Los Angeles,
Los Angeles, CA 90095}
\email{\{pola,tabuada\}@ee.ucla.edu}
\urladdr{http://www.ee.ucla.edu/$\sim$pola}
\urladdr{http://www.ee.ucla.edu/$\sim$tabuada}

\maketitle

\section{Introduction}
In recent years we have witnessed the development of different symbolic techniques aimed at reducing the complexity of controller synthesis~\cite{SymbolicIssue}. These techniques are based on the idea that many states can be treated as equivalent, when synthesizing controllers, and can thus be replaced by a symbol. The models resulting from replacing equivalent states by symbols, termed symbolic models, are typically simpler than the original ones, in the sense that they have a lower number of states. In many cases, one can even construct symbolic models with a finite number of states which is especially useful since controller design problems can then be solved on the symbolic models by resorting to well established results in supervisory control \cite{RamWonI} and algorithmic game theory \cite{Zielonka,AVW03}.
The search for classes of systems admitting symbolic models goes back to the 1990's and was motivated by problems of verification of dynamical and hybrid systems. Alur and Dill showed in \cite{TheoryTA} that timed automata admit symbolic models; this result was then generalized in \cite{AlurHS,Nicollin93} to multirate automata and in \cite{HenzingerDecHybSys,Varaiya94} to rectangular automata. More complex continuous dynamics, but simpler discrete dynamics, were considered in \cite{AHLP00}, where it was shown that o--minimal hybrid systems also admit symbolic models. Symbolic models for control systems were only considered later and early results were reported in \cite{KASL,TMoor,FJL02,CW98}. More precise results appeared recently in~\cite{LTLControl,SMC} where it was shown that discrete--time controllable linear systems admit symbolic models. Most of these results are based on appropriately adapting the notion of bisimulation introduced by Milner \cite{Milner} and Park \cite{Park} to the context of continuous and hybrid systems. A different approach emerged recently through the work of 
\cite{ApproxBis,RupakFORMATS05,AB-TAC07,Paulo_HSCC07}, where an approximate version of bisimulation was considered.  While (exact) bisimulation requires that observations of the states are identical, the notion of approximate bisimulation relaxes this condition, by allowing observations to be close and within a desired precision. This more flexible notion of bisimulation allows the identification of more classes of systems, admitting symbolic models. Indeed, the work in \cite{SCL} showed that for every asymptotically stabilizable linear control system it is possible to construct a symbolic model, which is based on an approximate notion of simulation (one--sided version of approximate bisimulation). Extensions of the results in \cite{SCL}, from approximate simulation to approximate bisimulation can be found in \cite{Girard_HSCC07,pola-2007}. In particular \cite{pola-2007} showed that, for the class of (incrementally globally) asymptotically stable nonlinear control systems,
symbolic models exist which are approximate bisimulation equivalent to control systems, with a precision that can be chosen a priori, as a design parameter.\\
Systems considered in the above described literature were either purely dynamical (e.g. \cite{AHLP00} and the references therein) or control systems (e.g. \cite{LTLControl, SCL,SMC,Girard_HSCC07,pola-2007}) not affected by exogenous disturbances. However, in many realistic situations, physical processes are characterized by a certain degree of uncertainty which is often modeled by additional exogenous disturbance inputs. 

\emph{The main contribution of this paper is to show that incrementally globally asymptotically stable control systems affected by exogenous inputs do admit symbolic models.}

The presence of disturbances requires us to replace the notion of approximate bisimulation used in \cite{AB-TAC07,pola-2007} with the notion of alternating approximate bisimulation, inspired by Alur and coworkers' alternating bisimulation \cite{Alternating}. To the best of the authors knowledge, alternating approximate bisimulation was never used before in the context of control systems. This novel notion of bisimulation is a critical ingredient of our results since, as illustrated in Section \ref{Sec32} through a simple example, approximate bisimulation fails to distinguish between the different role played by control inputs and disturbance inputs. Consequently, control strategies synthesized on symbolic models based on notions of bisimulation and approximate bisimulation cannot be transfered to the original models in a way which is robust with respect to disturbance inputs. Alternating approximate bisimulation solves this problem by guaranteeing that control strategies synthesized on symbolic models, based on alternating approximate bisimulations, can be readily transferred to the original model, independently of the particular evolution of the disturbance inputs.
In addition to show existence of symbolic models for a fairly general class of nonlinear control systems we also show that for linear control systems, symbolic models can be easily constructed by leveraging existing results on approximation of reachable sets (see e.g. \cite{Varaiya:98,Girard_HSCC05,Elipsoid,Tomlin_Reach} and the references therein).
Since control systems with disturbances can be thought of as arenas for differential games \cite{DifferentialGames}, our results also provide an alternative approach to the study of differential games by means of tools developed in computer science (see e.g. \cite{Zielonka,AVW03}).\\
Similar ideas to the ones of this paper have been recently explored in \cite{PolaTabuadaCDC07b} for the class of linear control systems with disturbances. A detailed discussion on relationships between results of the present paper and the ones in \cite{PolaTabuadaCDC07b}, can be found in the last section of this paper. A comparison with the work in \cite{BisimSchaft}, where systems with disturbances are also considered, appears in the last section of this paper.\\
This paper is organized as follows. Section \ref{sec2} introduces the class of control systems that we consider and some stability notions that will be used in the subsequent developments. Section \ref{sec3} introduces the notion of alternating transition systems that we use as an abstract representation of control systems and the notion of alternating approximate bisimulation upon which our results rely. Section \ref{sec4} is devoted to show existence of symbolic models for incrementally globally asymptotically stable nonlinear control systems. In Section \ref{sec5} we specialize the results of Section \ref{sec4} to the class of linear control systems and illustrate them in Section \ref{sec6}. Finally, some concluding remarks are offered in Section \ref{sec7}.

\section{Control systems and stability notions\label{sec2}}

\subsection{Notation}

The symbols $\mathbb{Z}$, $\mathbb{N}$, $\mathbb{R}$, $\mathbb{R}^{+}$ and $\mathbb{R}_{0}^{+}$ denote the set of integers, positive integers, reals, positive and nonnegative reals, respectively. The identity map on a set $A$ is denoted by $1_{A}$. 
Given two sets $A$ and $B$, if $A$ is a subset of $B$
we denote by \mbox{$\imath_{A}:A\hookrightarrow B$} or simply by $\imath$ the natural
inclusion map taking any $a\in A$ to \mbox{$\imath (a) = a \in B$}. Given a function $f:A\rightarrow B$ the symbol $f(A)$ denotes
the image of $A$ through $f$, i.e. $f(A):=\{b\in B:\exists a\in A$ s.t.
$b=f(a)\}$; if $C\subseteq A$, $f|_{C}:C\rightarrow B$ denotes the restriction of $f$ to $C$, so that $f|_{C}(c)=f(c)$ for any $c\in C$. We identify a relation
$R\subseteq A\times B$ with the map \mbox{$R:A\rightarrow2^{B}$} defined by
$b\in R(a)$ if and only if \mbox{$(a,b)\in R$}. Given a relation $R\subseteq A\times B$,
$R^{-1}$ denotes the inverse relation of $R$, i.e. \mbox{$R^{-1}:=\{(b,a)\in B\times A:(a,b)\in R\}$}.
Given a vector $x\in\mathbb{R}^{n}$ we denote by $x'$ the transpose of $x$ and by 
$x_{i}$ the $i$--th element of $x$; furthermore $\Vert x\Vert$ denotes the infinity norm of $x$; we recall that 
\mbox{$\Vert x\Vert:=\max\{|x_1|,|x_2|,...,|x_n|\}$}, where $|x_i|$ is the absolute value of $x_i$. Given a matrix $M$, the symbol $\Vert M\Vert$ denote the infinity norm of $M$; if \mbox{$M\in\mathbb{R}^{n\times m}$}, we
recall that $\Vert M\Vert:=\max_{1\leq i\leq m}%
{\textstyle\sum_{j=1}^{n}}
|a_{ij}|$. The symbol $conv(x^{1},x^{2},...,x^{m})$ denotes the
convex hull of vectors $x^{1},x^{2},...,x^{m}\in\mathbb{R}^{n}$. A bounded set
of the form $conv(x^{1},x^{2},...,x^{m})$ is called a polytope. Given a set $A\subseteq\mathbb{R}^{n}$, the symbol $\overline{A}$ denotes the topological closure of $A$. The symbol $\mathcal{B}_{\varepsilon}(x)$ denotes the closed ball centered at
\mbox{$x\in{\mathbb{R}}^{n}$} with radius $\varepsilon\in\mathbb{R}^{+}_{0}$, i.e. \mbox{$\mathcal{B}_{\varepsilon}(x)=\{y\in{\mathbb{R}}^{n} \,:\,\Vert x-y\Vert\leq\varepsilon\}$}. For any $A\subseteq
\mathbb{R}^{n}$ and \mbox{$\mu\in{\mathbb{R}}$} define \mbox{$[A]_{\mu}:=\{a\in A\,\,|$} \mbox{$a_{i}%
=k_{i}\mu,$}\mbox{$\,\,\,k_{i}\in\mathbb{Z}$} \mbox{$i=1,...,n\}$}. By geometrical considerations on the infinity norm, for any $\mu\in \mathbb{R}^{+}$ and $\lambda\geq\mu/2$ the collection of sets $\{\mathcal{B}_{\lambda}(q)\}_{q\in[\mathbb{R}^{n}]_{\mu}}$ is a covering of $\mathbb{R}^{n}$, i.e.
$\mathbb{R}^{n}\subseteq{\textstyle\bigcup\nolimits_{q\in\lbrack\mathbb{R}^{n}]_{\mu}}}\mathcal{B}_{\lambda}(q)$; conversely for any \mbox{$\lambda<\mu/2$}, \mbox{$\mathbb{R}^{n}\nsubseteq{\textstyle\bigcup\nolimits_{q\in\lbrack\mathbb{R}^{n}]_{\mu}}}\mathcal{B}_{\lambda}(q)$}. 
Given a measurable function \mbox{$f:\mathbb{R}_{0}^{+}\rightarrow\mathbb{R}$}, the
\mbox{(essential)} supremum of $f$ is denoted by $\Vert f\Vert_{\infty}$; we recall that $\Vert f\Vert_{\infty}:=(ess)sup\{\Vert f(t)\Vert,$ $t\geq0\}$; 
$f$ is essentially bounded if $\Vert f \Vert_{\infty} < \infty$. 
For a given time $\tau\in\mathbb{R}^{+}$, define $f_{\tau}$ so that
$f_{\tau}(t)=f(t)$, for any $t\in [0,\tau)$, and $f(t)=0$ elsewhere; 
$f$ is said to be locally essentially bounded if for any $\tau\in\mathbb{R}^{+}$,
$f_{\tau}$ is essentially bounded. A function \mbox{$f:\mathbb{R}^{n}\rightarrow \mathbb{R}$} is said to be radially unbounded if $f(x)\rightarrow \infty$, as $\Vert x\Vert\rightarrow \infty$. A continuous function $\gamma:\mathbb{R}_{0}^{+}%
\rightarrow\mathbb{R}_{0}^{+}$, is said to belong to class $\mathcal{K}$ if it
is strictly increasing and \mbox{$\gamma(0)=0$}; $\gamma$ is said to belong to class
$\mathcal{K}_{\infty}$ if \mbox{$\gamma\in\mathcal{K}$} and $\gamma(r)\rightarrow
\infty$, as $r\rightarrow\infty$. A continuous function \mbox{$\beta:\mathbb{R}_{0}^{+}\times\mathbb{R}_{0}^{+}\rightarrow\mathbb{R}_{0}^{+}$} is said to
belong to class $\mathcal{KL}$ if for each fixed $s$, the map $\beta(r,s)$
belongs to class $\mathcal{K}_{\infty}$ with respect to $r$ and, for each
fixed $r$, the map $\beta(r,s)$ is decreasing with respect to $s$ and
$\beta(r,s)\rightarrow0$, as \mbox{$s\rightarrow\infty$}. 
Given a metric space $(X,\mathbf{d})$, we denote by $\mathbf{d}_{h}$ the Hausdorff pseudo--metric induced by $\mathbf{d}$ on $2^{X}$; we recall that for any $X_{1},X_{2}\subseteq X$:
\[
\mathbf{d}_{h}(X_{1},X_{2}):=\max\{\vec{\mathbf{d}}_{h}(X_{1},X_{2}),\vec{\mathbf{d}}_{h}(X_{2},X_{1})\},
\]
where:
\[
\vec{\mathbf{d}}_{h}(X_{1},X_{2})=\sup_{x_{1}\in X_{1}}\inf_{x_{2}\in X_{2}}\mathbf{d}(x_{1},x_{2}),
\]
is the directed Hausdorff pseudo--metric. We recall that the Hausdorff pseudo--metric $\mathbf{d}_{h}$ satisfies the following properties for any $X_{1},X_{2},X_{3}\subseteq X$: (i) $X_{1}=X_{2}$ implies  $\mathbf{d}_{h}(X_{1},X_{2})=0$; (ii) $\mathbf{d}_{h}(X_{1},X_{2})=\mathbf{d}_{h}(X_{2},X_{1})$; (iii) $\mathbf{d}_{h}(X_{1},X_{3})\leq\mathbf{d}_{h}(X_{1},X_{2})+\mathbf{d}_{h}(X_{2},X_{3})$.

\subsection{Control Systems}

The class of systems that we consider in this paper is formalized in
the following definition.

\begin{definition}
\label{Def_control_sys}A \textit{control system} is a quadruple:
\[
\Sigma=(\mathbb{R}^{n},W,\mathcal{W},f),
\]
where:

\begin{itemize}
\item $\mathbb{R}^{n}$ is the state space;

\item $W=U\times V$ is the input space, where:

\begin{itemize}
\item $U\subseteq\mathbb{R}^{m}$ is the control input space;

\item $V\subseteq\mathbb{R}^{s}$ is the disturbance input space;
\end{itemize}

\item $\mathcal{W}=\mathcal{U}\times\mathcal{V}$ is a subset of the set of all
measurable and locally essentially bounded functions of time from intervals of the form $]a,b[\subseteq\mathbb{R}$ to $W$ with $a<0$ and $b>0$;

\item $f:\mathbb{R}^{n}\times W\rightarrow\mathbb{R}^{n}$ is a continuous map
satisfying the following Lipschitz assumption: for every compact set
$K\subset\mathbb{R}^{n}$, there exists a constant $L>0$ such that
\[
\Vert f(x,w)-f(y,w)\Vert\leq L\Vert x-y\Vert, 
\]
for all $x,y\in K$ and all $w\in W$. 
\end{itemize}
An absolutely continuous curve $\mathbf{x}:]a,b[\rightarrow\mathbb{R}^{n}$ is said
to be a \textit{trajectory} of $\Sigma$ if there exists $\mathbf{w}%
\in\mathcal{W}$ satisfying:
\[
\dot{\mathbf{x}}(t)=f(\mathbf{x}(t),\mathbf{w}(t)),
\]
for almost all $t\in$ $]a,b[$. 
\end{definition}
\bigskip

Although we have defined trajectories over open domains, we shall refer to
trajectories $\mathbf{x:}[0,\tau]\rightarrow\mathbb{R}^{n}$ defined on closed
domains $[0,\tau],$ $\tau\in\mathbb{R}^{+}$ with the understanding of the
existence of a trajectory $\mathbf{z}:]a,b[\rightarrow\mathbb{R}^{n}$
such that $\mathbf{x}=\mathbf{z}|_{[0,\tau]}$. We will also write
$\mathbf{x}(\tau,x,\mathbf{w})$ to denote the point reached at time $\tau \in ]a,b[$
under the input $\mathbf{w}$ from initial condition $x$; this point is
uniquely determined, since the assumptions on $f$ ensure existence and
uniqueness of trajectories. \\
In some of the subsequent developments we assume that control systems are forward complete. 
We recall that a control system $\Sigma$ is \textit{forward complete} if every trajectory is defined on an interval of the
form $]a,\infty\lbrack$. The following result completely characterizes forward completeness.

\begin{theorem} \cite{fc-theorem}
Consider a control system $\Sigma=(\mathbb{R}^{n},W,\mathcal{W},f)$ and suppose that $W$ is compact.
Then $\Sigma$ is forward complete if and only if there exists a radially unbounded smooth function 
\mbox{$\mathbf{V}:\mathbb{R}^{n}\rightarrow\mathbb{R}_{0}^{+}$} such that for any $x \in \mathbb{R}^{n}$ and for any $w \in W$ the following exponential growth condition is verified:
\[
\frac{\partial \mathbf{V}}{\partial x}f(x,w) \leq \mathbf{V}(x).
\]
\end{theorem}
Simpler, but only sufficient, conditions for forward completeness are also available in the literature. These include linear growth or compact support of the vector field (see e.g. \cite{OCT}). Whenever we need to distinguish between a control input value $u$ and a disturbance input value $v$ in $(u,v)\in W$ we
slightly abuse notation by writing $f(x,u,v)$ instead of $f(x,(u,v))$. Analogously, whenever we need to distinguish between $\mathbf{u}$ and $\mathbf{v}$ in an input signal
$(\mathbf{u},\mathbf{v})\in \mathcal{W}$, we write $\mathbf{x}%
(\tau,x,\mathbf{u},\mathbf{v})$ instead of $\mathbf{x}(\tau
,x,(\mathbf{u},\mathbf{v}))$.  

\subsection{Stability notions}

The results presented in this paper will assume certain stability assumptions
on the control systems. We briefly recall those notions 
and results that will be used in this paper.

\begin{definition}
\label{Def_IRAS} \cite{IncrementalS} A control system $\Sigma$ is said to be
\textit{incrementally globally asymptotically stable} ($\delta$--GAS) if it is
forward complete and there exist a $\mathcal{KL}$ function $\beta$ such that
for any $t\in \mathbb{R}^{+}_{0}$, any $x_{1},x_{2}\in\mathbb{R}^{n}$ and any input signal
$\mathbf{w}\in\mathcal{W}$ the following condition is satisfied:%
\begin{equation}
\left\Vert \mathbf{x}(t,x_{1},\mathbf{w})-\mathbf{x}(t,x_{2},\mathbf{w}%
)\right\Vert \leq\beta(\left\Vert x_{1}-x_{2}\right\Vert ,t).
\label{deltaUGAS}%
\end{equation}
\end{definition}

The above definition can be thought of as an incremental version of the classical notion of
global asymptotic stability (GAS) \cite{Khalil}. Furthermore when $f$ satisfies $f(0,0)=0$, $\delta$--GAS implies GAS of $\Sigma$ with $W=\{0\}$, by just comparing a trajectory of $\Sigma$ with any initial condition $x\in\mathbb{R}^{n}$ and identically null input $\mathbf{w}(t)=0$, $t\in\mathbb{R}^{+}_{0}$, with the null trajectory $\mathbf{x}(t)=0,$ $t\in\mathbb{R}^{+}_{0}$.
In general, it is difficult to check directly inequality (\ref{deltaUGAS}). However, $\delta$--GAS can be characterized by dissipation inequalities.

\begin{definition}
\cite{IncrementalS} Given a control system $\Sigma$, a smooth function:
\[
\mathbf{V}:\mathbb{R}^{n}\times\mathbb{R}^{n}\rightarrow\mathbb{R}_{0}^{+},
\]
is called a $\delta$--GAS\textit{ Lyapunov function} for $\Sigma$, if there exist
$\mathcal{K}_{\infty}$ functions $\alpha_{1}$, $\alpha_{2}$ and $\rho$\ such that:

\begin{itemize}
\item[(i)] for any $(x_{1},x_{2})\in\mathbb{R}^{n}\times\mathbb{R}^{n}$%
\[
\alpha_{1}(\left\Vert x_{1}-x_{2}\right\Vert )\leq \mathbf{V}(x_{1},x_{2})\leq
\alpha_{2}(\left\Vert x_{1}-x_{2}\right\Vert )\text{;}%
\]

\item[(ii)] for any $x_{1},x_{2}\in\mathbb{R}^{n}$ and any $w\in W$%
\[
\frac{\partial \mathbf{V}}{\partial x_{1}}f(x_{1},w)+\frac{\partial \mathbf{V}}{\partial x_{2}%
}f(x_{2},w)<-\rho(\left\Vert x_{1}-x_{2}\right\Vert )\text{.}%
\]

\end{itemize}
\end{definition}
The following result completely characterizes $\delta$--GAS of a control system in terms of existence of Lyapunov functions.

\begin{theorem}
\cite{IncrementalS} Consider a forward complete control system $\Sigma=(\mathbb{R}%
^{n},W,$ $\mathcal{W},f)$ and suppose that $W$ is a compact subset of $\mathbb{R}%
^{m} \times \mathbb{R}^{s}$. Then $\Sigma$ is $\delta$--GAS if and only if it admits a $\delta$--GAS Lyapunov function.
\end{theorem}

\section{Symbolic models and approximate equivalence notions}\label{sec3}

\subsection{Alternating transition systems}

In this paper we will use the class of alternating transition systems as abstract models of control systems.

\begin{definition}
\label{Def_TS}An (alternating) transition system is a tuple:
\[
T=(Q,L,\rTo,O,H),
\]
consisting of:

\begin{itemize}
\item A set of states $Q$;

\item A set of labels $L=A\times B$, where:

\begin{itemize}
\item $A$ is the set of control labels;

\item $B$ is the set of disturbance labels;

\end{itemize}

\item A transition relation $\rTo\subseteq Q\times L\times Q$;

\item An output set $O$;

\item An output function $H:Q\rightarrow O$.
\end{itemize}

A transition system $T$ is said to be:

\begin{itemize}
\item \textit{metric}, if the output set $O$ is equipped with a metric
$\mathbf{d}:O\times O\rightarrow\mathbb{R}_{0}^{+}$;

\item \textit{countable}, if $Q$ and $L$ are countable sets;

\item \textit{finite}, if $Q$ and $L$ are finite sets.
\end{itemize}
\end{definition}
\bigskip
We will follow standard practice and denote by $q \rTo^{a,b} p$, a transition from $q$ to $p$
labeled by $a,b$. Transition
systems capture dynamics through the transition relation. For any states
$q,p\in Q$, $q \rTo^{a,b} p$ simply means that it is
possible to evolve or jump from state $q$ to state $p$ under the action
labeled by $a$ and $b$. A transition system can be represented as a graph where circles represent states and arrows represent transitions (see e.g. Figure \ref{fig11}). A transition system as in Definition \ref{Def_TS}, can
be thought of as an arena for a $2$--players game, where the protagonist acts by choosing
control labels and the antagonist acts by choosing disturbance labels. 
We will use transition systems as an abstract representation of control systems. 
There are several different ways in which we can transform control
systems into transition systems. We now describe one of these which has the
property of capturing all the information contained in a control system
$\Sigma$. Given $\Sigma=(\mathbb{R}^{n},U \times V,\mathcal{U} \times \mathcal{V},f)$ define the
transition system:
\[
T(\Sigma):=(Q,L,\rTo,O,H),
\]
where:

\begin{itemize}
\item $Q=\mathbb{R}^{n}$;

\item $L=A \times B$ where $A=\mathcal{U}$ and $B= \mathcal{V}$;

\item $q\rTo^{\mathbf{u,v}}p$ if a trajectory
$\mathbf{x}:[0,\tau]\rightarrow\mathbb{R}^{n}$ of $\Sigma$ exists so that
$\mathbf{x}(\tau,q,\mathbf{u,v})=p$ for some $\tau\in\mathbb{R}^{+}$;

\item $O=\mathbb{R}^{n}$;

\item $H=1_{\mathbb{R}^{n}}$.
\end{itemize}
\bigskip
In the subsequent developments we will work with a sub--transition system of $T(\Sigma)$
obtained by selecting those transitions from $T(\Sigma)$ describing
trajectories of duration $\tau$ for some chosen $\tau\in\mathbb{R}^{+}$. This
can be seen as a time discretization or sampling process.

\pagebreak

\begin{definition}
Given a control system $\Sigma=(\mathbb{R}^{n},U\times V,\mathcal{U}%
\times\mathcal{V},f)$ and a parameter $\tau\in\mathbb{R}^{+}$ define the  
transition system:%
\[
T_{\tau}(\Sigma):=(Q_{\tau},L_{\tau},\rTo_{\tau},O_{\tau},H_{\tau}),
\]
where:

\begin{itemize}
\item $Q_{\tau}=\mathbb{R}^{n}$;

\item $L_{\tau}=A_{\tau}\times B_{\tau}$ where:
\begin{eqnarray}
&&A_{\tau}=\{\mathbf{u}\in\mathcal{U}\,|\,\text{the domain of $\mathbf{u}$ is $[0,\tau]$}\},\notag\\
&&B_{\tau}=\{\mathbf{v}\in\mathcal{V}\,|\,\text{the domain of $\mathbf{v}$ is $[0,\tau]$}\};\notag
\end{eqnarray}
 
\item $q\rTo^{\mathbf{u,v}}_{\tau} p$ if a
trajectory $\mathbf{x}:[0,\tau]\rightarrow\mathbb{R}^{n}$ of $\Sigma$ exists
so that $\mathbf{x}(\tau,q,\mathbf{u,v})=p$;

\item $O_{\tau}=\mathbb{R}^{n}$;

\item $H_{\tau}=1_{\mathbb{R}^{n}}$.
\end{itemize}
\end{definition}
\bigskip
Note that $T_{\tau}(\Sigma)$ is a metric transition system when we regard
$O_{\tau}=\mathbb{R}^{n}$ as being equipped with the metric $\mathbf{d}(p,q)=\left\Vert
p-q\right\Vert $. 

\subsection{Alternating and approximate bisimulations}\label{Sec32}
In this section we introduce a notion of approximate equivalence upon which
all the results in this paper rely. The notion that we consider, is the one of bisimulation equivalence
\cite{Milner,Park}. Bisimulation relations are standard mechanisms to relate
the properties of transition systems\ \cite{ModelChecking}. Intuitively, a
bisimulation relation between a pair of transition systems $T_{1}$ and $T_{2}$
is a relation between the corresponding sets of states explaining how a sequence
of transitions $r_{1}$ of $T_{1}$ can be transformed into a sequence of
transitions $r_{2}$ of $T_{2}$ and vice versa. While typical bisimulation
relations require that $r_{1}$ and $r_{2}$ are observationally
indistinguishable, that is $H_{1}(r_{1})=H_{2}(r_{2})$, we shall relax this
by requiring $H_{1}(r_{1})$ to simply be close to $H_{2}(r_{2})$, where
closeness is measured with respect to the metric on the output set. 
The following definition has been introduced in \cite{AB-TAC07} and in a slightly different formulation in \cite{Paulo_HSCC07}.
\begin{definition}
\label{ASR_bis}
Given two transition systems $T_{1}=(Q_{1},L_{1},\rTo_{1},O,H_{1})$ and
\mbox{$T_{2}=(Q_{2},L_{2},\rTo_{2},O,H_{2})$} with the same output set and metric $\mathbf{d}$, and given a precision $\varepsilon \in \mathbb{R}_{0}^{+}$, a relation 
\[
R\subseteq Q_{1}\times Q_{2},
\]
is said to be an $\varepsilon
$--approximate bisimulation relation between $T_{1}$ and $T_{2}$, if for
any $(q_{1},q_{2})\in R$:
\begin{itemize}
\item[(i)]$\mathbf{d}(H_{1}(q_{1}),H_{2}(q_{2}))\leq\varepsilon$;
\item[(ii)]$q_{1}\rTo^{l_{1}}_{1} p_{1}$ implies 
existence of $q_{2}\rTo^{l_{2}}_{2} p_{2}$ such that
$(p_{1},p_{2})\in R$;
\item[(iii)]$q_{2}\rTo^{l_{2}}_{2} p_{2}$ implies 
existence of $q_{1}\rTo^{l_{1}}_{1} p_{1}$ such that
$(p_{1},p_{2})\in R$.
\end{itemize}
Moreover $T_{1}$ is $\varepsilon$--approximately bisimilar to $T_{2}$ if there exists
an \mbox{$\varepsilon$--approximate} bisimulation relation $R$\ between $T_{1}$ and $T_{2}$ such that
\mbox{$R(Q_{1})=Q_{2}$} and $R^{-1}(Q_{2})=Q_{1}$.
\end{definition}
\bigskip

Note that when $\varepsilon=0$, the notion of $\varepsilon$--approximate bisimulation relation is substantially equivalent to the classical notion of Milner \cite{Milner} and Park \cite{Park}. 
The work in \cite{pola-2007} showed existence of symbolic models that are approximately bisimilar to $\delta$--GAS control systems (with no disturbance). However, the notion in Definition \ref{ASR_bis} employed in \cite{pola-2007}, does not capture the different role of control and disturbance inputs in control systems. The following example shows that approximate bisimulation (in the sense of Definition \ref{ASR_bis}) cannot be used for control design of systems affected by disturbances.
\begin{example}
\label{exa}
Consider the following control system:
\begin{equation}
\Sigma=(\mathbb{R},U\times V,\mathcal{U}\times \mathcal{V},f), 
\label{example}
\end{equation}
where \mbox{$U=[1,2]\subset \mathbb{R}$}, \mbox{$V=[0.4,1]\subset \mathbb{R}$}, $\mathcal{U}\times \mathcal{V}$ is the class of all measurable and locally essentially bounded functions taking values in $U\times V$, and $f:\mathbb{R}\times U\times V\rightarrow\mathbb{R}$ is defined by $f(x,u,v)=-2x+uv$. We work in the compact state space\footnote{The set $X$ is invariant for the control system $\Sigma$, i.e. $\mathbf{x}(t,x,\mathbf{u},\mathbf{v})\in X$, for any $x\in X$, any $(\mathbf{u},\mathbf{v})\in \mathcal{U}\times\mathcal{V}$, and any time $t\in\mathbb{R}^{+}_{0}$. Indeed it is easy to see that for any state $x$ in the boundary of $X$ and for any input $(u,v)\in U \times V$, $f$ points in $X$, i.e. $f(x,u,v)\in X$.} $X=[0,2]$.
Consider the transition system:
%
\begin{equation}
T=(Q,L,\rTo,O,H),
\label{T2example}
\end{equation}
where:
\begin{itemize}
\item $Q=\{q_{1},q_{2},q_{3}\}$;

\item $L=\{l_{1},l_{2},l_{3}\}$;

\item $q\rTo^{l}p$ is depicted in Figure \ref{fig11};

\item $O=\mathbb{R}$;

\item $H:O\rightarrow \mathbb{R}$ is defined by $H(q_{1})=0$, $H(q_{2})=1$, and $H(q_{3})=2$.
\end{itemize}
Given the desired precision $\varepsilon=0.6$ and $\tau=1$, by using the results in \cite{pola-2007}, it is possible to show that the relation $R\subset Q_{\tau}\times Q$ defined by:
\begin{equation}
R=R_{1}\times\{q_{1}\}\cup R_{2}\times\{q_{2}\}\cup R_{3}\times\{q_{3}\},
\label{RR_example}
\end{equation}
where $R_{1}=[0,0.6]$, $R_{2}=[0.4,1.6]$ and $R_{3}=[1.4,2]$, is a $0.6$--approximate bisimulation relation between $T_{\tau}(\Sigma)$ and $T$. Furthermore, since $R(Q_{\tau})=Q$ and $R^{-1}(Q)=Q_{\tau}$ transition systems $T_{\tau}(\Sigma)$ and $T$ are $0.6$--approximately bisimilar\footnote{Transition system $T$ coincides with transition system $T_{\tau,\eta,\mu}(\Sigma)$ as defined in (4.4) of \cite{pola-2007}, with $\tau=1$, $\eta=1$ and $\mu=0.01$. Theorem 4.2 of \cite{pola-2007} guarantees that $T$ is $0.6$--approximately bisimilar to transition system $T_{\tau}(\Sigma)$ with $\tau=1$. Notice that condition (4.5) of \cite{pola-2007} boils down, in this case, to $e^{-2\tau}\varepsilon+\mu+\eta/2\leq\varepsilon$, which is indeed satisfied.}.
Suppose now that the goal is to find a control strategy on $T$ such that, starting from state $q_{1}$ it is possible to reach the set $\{q_{2},q_{3}\}$ in one step.
By Figure \ref{fig11}, $q_{1}\rTo^{l_{2}}q_{2}$ and $q_{1}\rTo^{l_{3}}q_{3}$ and hence both labels $l_{2}$ and $l_{3}$ solve that problem. Since $(0,q_{1})\in R$, the notion of approximate bisimulation (see condition (iii) of Definition \ref{ASR_bis}) guarantees that starting from $0\in R_{1}$ there exists a pair of labels $(a_{2},b_{2}),(a_{3},b_{3})\in A_{\tau} \times B_{\tau}$ so that: 
\begin{eqnarray}
0\rTo^{a_{2},b_{2}}_{\tau}x_{2}\in R_{2},
&
\hspace{5mm}
&
0\rTo^{a_{3},b_{3}}_{\tau}x_{3}\in R_{3},
\label{zaq}
\end{eqnarray}
in transition system $T_{\tau}(\Sigma)$. Indeed, by choosing constant curves $(a_{2}(t),b_{2}(t))=(1,1)$ and $(a_{3}(t),b_{3}(t))=(2,1)$, $t\in [0,1]$ we have:
\begin{eqnarray}
0\rTo^{a_{2},b_{2}}_{\tau}0.865\in R_{2},
&
\hspace{5mm}
&
0\rTo^{a_{3},b_{3}}_{\tau}1.73\in R_{3}.\notag
\end{eqnarray}
However, if the constant disturbance label $b(t)=0.4,t\in [0,1]$ occurs instead of \mbox{$b_{2}=b_{3}$}, we obtain:
\begin{eqnarray}
0\rTo^{a_{2},b}_{\tau}0.346\in R_{1},
&
\hspace{5mm}
&
0\rTo^{a_{3},b}_{\tau}0.692\in R_{2},
\end{eqnarray}
thus showing that the control strategy in (\ref{zaq}) does not produce the desired result on the transition system $T_{\tau}(\Sigma)$.
Although $T$ is not adequate to solve this problem, a solution does exist.
Since $0\rTo^{a_{3},b}_{\tau}0.692\in R_{2}$ and the set $X$ is invariant for $\Sigma$, it is easy to see that for any $\hat{b}\in B_{\tau}$, $0\rTo^{a_{3},\hat{b}}_{\tau}x$ with $x\geq 0.692$ and hence $x\in R_{2}\cup R_{3}$. Therefore, control label $a_{3}$ guarantees that state $0\in R_{1}$ reaches $R_{2}\cup R_{3}$, robustly with respect to the disturbance labels action, whereas control label $a_{2}$ does not. We stress that this different feature of control labels $a_{2}$ and $a_{3}$ is not captured by the notion of approximate bisimulation in Definition \ref{ASR_bis}.

\begin{figure}[ptb]
\begin{center}
\begin{tikzpicture}[->,shorten >=1pt,auto,node distance=3cm,semithick,inner sep=1.5pt,bend angle=45]
\node[state] (q_1) {$q_{1}$};
\node[state] (q_2) [ right of=q_1] {$q_{2}$};
\node[state] (q_3) [ right of=q_2] {$q_{3}$};
\path[->] (q_1) edge  [bend left]  node[below] {$l_{2}$} (q_2)
edge [bend left] node [below] {$l_{3}$} (q_3)
edge [loop left] node {$l_{1}$} ()
(q_2) edge [bend left] node [below] {$l_{3}$} (q_3)
edge [loop left] node {$l_{1},l_{2}$} ()
(q_3) edge  [bend left]  node [above] {$l_{1},l_{2}$} (q_2)
edge [loop right] node {$l_{3}$} ();
\end{tikzpicture}
\end{center}
\caption{Transition system $T$ defined in (\ref{T2example}) associated with control system $\Sigma$ defined in (\ref{example}).}
\label{fig11}%
\end{figure}
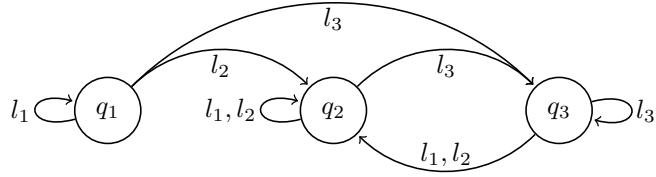
\label{qwert}
\end{example}
\bigskip

The above example motivates us to propose the following definition that combines the notions of \cite{AB-TAC07} and \cite{Paulo_HSCC07}, with the notion of alternating bisimulation, introduced by Alur and coworkers in \cite{Alternating}.
\begin{definition}
\label{Alt_ASR}Given two transition system $T_{1}=(Q_{1},A_{1}\times B_{1},\rTo_{1},O,H_{1})$ and $T_{2}=(Q_{2},A_{2}\times B_{2},\rTo_{2},O,H_{2})$ with the same observation set and the same metric $\mathbf{d}$
and given a precision $\varepsilon\in \mathbb{R}^{+}_{0}$, a relation 
\[
R\subseteq Q_{1}\times Q_{2},
\]
is said to be an alternating $\varepsilon$--approximate ($A\varepsilon A$)
bisimulation relation between $T_{1}$ and $T_{2}$ if for any $(q_{1}%
,q_{2})\in R$:
\begin{itemize}
\item[(i)]$\mathbf{d}(H_{1}(q_{1}),H_{2}(q_{2}))\leq\varepsilon$;
\item[(ii)] $\forall a_{1}$ $\in A_{1}$ $\exists a_{2}\in A_{2}$ $\forall b_{2}\in B_{2}$ $\exists b_{1}\in B_{1}$ such that:
\begin{eqnarray}
q_{1} \rTo^{a_{1},b_{1}}_{1} p_{1},
&
\hspace{5mm}
&
q_{2} \rTo^{a_{2},b_{2}}_{2} p_{2},
\label{qqq}
\end{eqnarray}
with $(p_{1},p_{2})\in R$;
\item[(iii)]$\forall a_{2}\in A_{2}$ $\exists a_{1}\in A_{1}$ $\forall b_{1}\in B_{1}$ $\exists b_{2}\in B_{2}$ such that: 
\begin{eqnarray}
q_{1} \rTo^{a_{1},b_{1}}_{1} p_{1},
&
\hspace{5mm}
&
q_{2} \rTo^{a_{2},b_{2}}_{2} p_{2},
\label{qqq1}
\end{eqnarray}
with $(p_{1},p_{2})\in R$.
\end{itemize}

Moreover, $T_{1}$ is said to be $A\varepsilon A$ bisimilar to $T_{2}$ if there exists
an $A\varepsilon A$ bisimulation relation $R$ between $T_{1}$
and $T_{2}$ such that $R(Q_{1})=Q_{2}$ and $R^{-1}(Q_{2})=Q_{1}$.
\end{definition}
\bigskip

It is easy to see that Definition \ref{ASR_bis} can be recovered as a special case of Definition \ref{Alt_ASR}, when the cardinality of each of the sets $B_{1}$ and $B_{2}$ in transition systems $T_{1}$ and $T_{2}$ is one. Moreover when $\varepsilon=0$, the notion of bisimulation in Definition \ref{Alt_ASR} coincides with the 2--players version of the definition proposed in \cite{Alternating}.\\
Definition \ref{Alt_ASR} captures the different role played by control and disturbance labels in the transition systems involved, whereas Definition \ref{ASR_bis} does not. In fact, by \cite{Alternating} it is possible to show that $A\varepsilon A$ bisimulation relations \textit{preserve control strategies} (see Lemma 1 in \cite{Alternating}) and hence prevent phenomena illustrated in Example \ref{qwert}. Indeed, conditions (ii) and (iii) of Definition \ref{Alt_ASR} require that the choice of control labels is made ``robustly'' with respect to the action of disturbance labels. For example, given any control label $a_{1}$ in $T_{1}$, condition (ii) requires existence of a control label $a_{2}$ in $T_{2}$ so that, for any possible action of the disturbance labels $b_{1}$ in $T_{1}$ and $b_{2}$ in $T_{2}$ respectively, the corresponding transitions $q_{1} \rTo^{a_{1},b_{1}}_{1} p_{1}$ and
$q_{2} \rTo^{a_{2},b_{2}}_{2} p_{2}$ match the relation $R$, i.e. $(p_{1},p_{2})\in R$.\\
A dual notion of $A\varepsilon A$ bisimulation relation can be given when we reverse the role of control and disturbance labels in conditions (ii) and (iii), as follows:
\bigskip
\begin{itemize}
\item[(ii')] 
$\forall b_{1}\in B_{1}$ $\exists b_{2}\in B_{2}$ $\forall a_{2}\in A_{2}$ $\exists a_{1}\in A_{1}$ such that (\ref{qqq}) holds, with $(p_{1},p_{2})\in
R$;
\item[(iii')] 
$\forall b_{2}\in B_{2}$ $\exists b_{1}\in B_{1}$ $\forall a_{1}\in A_{1}$ $\exists a_{2}\in A_{2}$ such that (\ref{qqq1}) holds, with $(p_{1},p_{2})\in
R$.
\end{itemize}
\bigskip
For later use, whenever we want to distinguish between the two notions of $A\varepsilon A$ bisimulation relation, we refer to the notion of Definition \ref{Alt_ASR} as $(A_{1},A_{2})$--$A\varepsilon A$ bisimulation relation and to the notion of Definition \ref{Alt_ASR}, where conditions (ii) and (iii) are replaced by conditions (ii') and (iii'), by $(B_{1},B_{2})$--$A\varepsilon A$ bisimulation relation. Furthermore, a bisimulation relation $R$ that satisfies conditions (i), (ii), (iii), (ii') and (iii') is called a $(A_{1},A_{2})$--$(B_{1},B_{2})$--$A\varepsilon A$ bisimulation relation. Consequently, if $R(Q_{1})=Q_{2}$ and $R^{-1}(Q_{2})=Q_{1}$, transition systems $T_{1}$ and $T_{2}$ are said to be $(A_{1},A_{2})$--$(B_{1},B_{2})$--$A\varepsilon A$ bisimilar.\\

\section{Existence of symbolic models\label{sec4}}

In this section we present the main result of this paper:
\begin{theorem}
\label{Th_main22}Consider a control system $\Sigma=(\mathbb{R}^{n},U\times
V,\mathcal{U\times V},f)$. If $\Sigma$ is \mbox{$\delta$--GAS} and $U\times V$ is compact, then for any desired precision $\varepsilon\in
\mathbb{R}^{+}$ there exist $\tau\in \mathbb{R}^{+}$ and a countable transition system $T$ that is $A\varepsilon A$ bisimilar to $T_{\tau}(\Sigma)$.
\end{theorem}
\bigskip

The above result is important because it shows existence of symbolic models for nonlinear control systems \textit{in presence of disturbances}, and therefore it provides a first step toward the construction of symbolic models with guaranteed approximation properties.
In fact, in the next section we show how to construct countable transition systems that are $A\varepsilon A$ bisimilar to \textit{linear} control systems.
\\
Theorem \ref{Th_main22} relies upon the $\delta$--GAS assumption on the control system considered. This condition is not far from also being necessary. 
The following counterexample shows that unstable control systems do not admit, in general, countable symbolic models.
\begin{example}
Consider a control system $\Sigma=(\mathbb{R},U \times V,\mathcal{U}\times \mathcal{V},f)$, where \mbox{$U\times V=\{0\}\times \{0\}$}, $\mathcal{U}\times \mathcal{V}=\{\mathbf{0}\}\times \{\mathbf{0}\}$, $\mathbf{0}$ is the identically null input and $f(x)=x$. The input space $U\times V$ is compact; furthermore $\Sigma$ is unstable and hence not $\delta$--GAS. Hence, $\Sigma$ satisfies all the conditions required in Theorem \ref{Th_main22}, except for $\delta$--GAS. We now show that there exists an $\varepsilon\in \mathbb{R}^{+}$ such that for any $\tau\in \mathbb{R}^{+}$ and any countable transition system $T$, transition systems $T_{\tau}(\Sigma)$ and $T$ are \textit{not} $A\varepsilon A$ bisimilar. 
Since $U\times V=\{0\}\times \{0\}$, the notions of bisimulation in Definition \ref{ASR_bis} and in Definition \ref{Alt_ASR} coincide, and therefore in the following we will work, for simplicity, with the one in Definition \ref{ASR_bis}. Pick any $\varepsilon\in \mathbb{R}^{+}$, any $\tau\in \mathbb{R}^{+}$ and any countable metric transition system:
\[
T=(Q,L,\rTo,\mathbb{R},H),
\]
with $H:Q\rightarrow \mathbb{R}$ and the same metric $\mathbf{d}(p,q)=\Vert p-q \Vert$ of $T_{\tau}(\Sigma)$. Consider any relation $R\subseteq Q_{\tau}\times Q$ satisfying conditions (i), (ii) and (iii) of Definition \ref{ASR_bis} and such that $R(Q_{\tau})=Q$ and $R^{-1}(Q)=Q_{\tau}$. We now show that such relation $R$ does not exist.
By countability of $T$, there exist $q_{0}\in Q$ and $x_{0},y_{0}\in Q_{\tau}=\mathbb{R}$ such that $x_{0}\neq y_{0}$, and $(x_{0},q_{0}),(y_{0},q_{0})\in R$. Set $x_{k}=e^{\tau k}x_{0}$, $y_{k}=e^{\tau k}y_{0}$, for any $k\in\mathbb{N}$. 
Since $x_{0}\neq y_{0}$, by selecting $\lambda\in\mathbb{R}^{+}$ such that $\Vert x_{0}-y_{0} \Vert>\lambda$, we have:
\begin{equation}
\Vert x_{k}-y_{k} \Vert = e^{\tau k}\Vert x_{0}-y_{0} \Vert > e^{\tau k}\lambda,\forall k\in\mathbb{N}.
\label{zexample}
\end{equation}
Choose $k' \in \mathbb{N}$ so that $e^{\tau k'} \lambda - \varepsilon>\varepsilon$. By conditions (ii) and (iii) in Definition \ref{ASR_bis} and since $R(Q_{\tau})=Q$ and $R^{-1}(Q)=Q_{\tau}$, there must exist $q_{k'}\in Q$ so that, $(x_{k'},q_{k'}),(y_{k'},q_{k'})\in R$. Since $(x_{k'},q_{k'})\in R$,
\begin{equation}
\Vert x_{k'}-H(q_{k'}) \Vert \leq \varepsilon.
\label{xexample}
\end{equation}
By combining inequalities (\ref{zexample}) and (\ref{xexample}) and by definition of $k'$, we obtain:
\begin{equation}
\Vert H(q_{k'})-y_{k'} \Vert \geq \Vert x_{k'} - y_{k'} \Vert -\Vert x_{k'}- H(q_{k'})\Vert > e^{\tau k'} \lambda - \varepsilon>\varepsilon.
\label{zzexample}
\end{equation}
Inequality (\ref{zzexample}) shows that the pair $(y_{k'},q_{k'})\in R$ does not satisfy condition (i) of Definition \ref{ASR_bis} and hence we conclude that a relation $R\subseteq Q_{\tau}\times Q$ satisfying conditions (i), (ii) and (iii) of Definition \ref{ASR_bis} and such that $R(Q_{\tau})=Q$ and $R^{-1}(Q)=Q_{\tau}$ does not exist. Thus transition systems $T_{\tau}(\Sigma)$ and $T$ are \textit{not} $A\varepsilon A$ bisimilar.
\end{example}
\bigskip

The last part of this section will be devoted to the proof of Theorem \ref{Th_main22}, which is based on three steps:
\begin{itemize}
\item[(1)] we first associate a suitable transition system $T_{\tau,\eta,\mu}(\Sigma)$ to a control system $\Sigma=(\mathbb{R},U\times V,\mathcal{U}\times \mathcal{V},f)$ (Definition \ref{SymbModNLCS} on page 13);
\item[(2)] we then prove, under a compactness assumption on $U\times V$, that transition system $T_{\tau,\eta,\mu}(\Sigma)$ is countable (Corollary \ref{Corollary1} on page 14);
\item[(3)] we finally prove, under the $\delta$--GAS assumption on $\Sigma$, that $T_{\tau,\eta,\mu}(\Sigma)$ is $A\varepsilon A$ bisimilar to $T_{\tau}(\Sigma)$ (Theorem \ref{Th_main2} on page 14).
\end{itemize}
\bigskip

\textbf{STEP 1.}
Given a control system $\Sigma$, any $\tau\in\mathbb{R}^{+}$, $\eta
\in\mathbb{R}^{+}$ and $\mu\in\mathbb{R}^{+}$ we will define the transition
system:
\begin{equation}
T_{\tau,\eta,\mu}(\Sigma):=(Q,L,\rTo,O,H). 
\label{T4}
\end{equation}
Parameters $\tau\in\mathbb{R}^{+},\eta\in\mathbb{R}^{+}$ and $\mu\in
\mathbb{R}^{+}$ in transition system $T_{\tau,\eta,\mu}(\Sigma)$ can be
thought of, respectively, as a sampling time, a state space and an input space
quantization. In order to define $T_{\tau,\eta,\mu}(\Sigma)$ we will extract a
countable set of states $Q$ from $Q_{\tau}$ and a countable set of labels $L$ from $L_{\tau}$, in such a way, that the resulting $T_{\tau,\eta,\mu}(\Sigma)$ is countable and $A\varepsilon A$ bisimilar to $T_{\tau}(\Sigma)$.

From now on, we denote by $\mathbf{d}_{h}$ the Hausdorff pseudo--metric induced by
the metric $\mathbf{d}$ of the observation space $O_{\tau}$ of $T_{\tau}(\Sigma)$. Furthermore, since the output function $H_{\tau}$ of $T_{\tau}(\Sigma)$ is the identity function, we write $\mathbf{d}(x,y)=\Vert x-y \Vert$, instead of 
\mbox{$\mathbf{d}(x,y)=\Vert H_{\tau}(x)-H_{\tau}(y) \Vert$}.\\
We start by showing that any subset of $\mathbb{R}^{n}$ can be arbitrarily
well approximated by a subset of the lattice $[\mathbb{R}^{n}]_{\eta}$, where
$\eta$ is the precision that we require on the approximation.

\begin{lemma}
\label{Lemma1}For any set $X\subseteq\mathbb{R}^{n}$ and any precision $\eta\in
\mathbb{R}^{+}$ there exists $P\subseteq\lbrack\mathbb{R}^{n}]_{\eta}$ such
that $\mathbf{d}_{h}(P,X)\leq\eta/2$.
\end{lemma}
\begin{proof}
By geometrical considerations on the infinity norm, $X\subseteq%
{\textstyle\bigcup\nolimits_{p\in\lbrack\mathbb{R}^{n}]_{\eta}}}
\mathcal{B}_{\eta/2}(p)$ and therefore for
any $x\in X$ there exists $p\in\lbrack\mathbb{R}^{n}]_{\eta}$ such that
$\mathbf{d}(x,p)=\Vert x-p\Vert
\leq\eta/2$. Denote by $\vartheta:X\rightarrow\lbrack\mathbb{R}^{n}]_{\eta}$ a
function that associates to any $x\in X$ a vector $p\in\lbrack\mathbb{R}%
^{n}]_{\eta}$ so that $\mathbf{d}(x,p)=\Vert x-p\Vert\leq\eta/2$ and set
$P=\vartheta(X)$. Notice that by construction, for any $p\in P$ there exists
$x\in X$ such that $\mathbf{d}(x,p)=\Vert x-p\Vert\leq\eta/2$ (choose $x$ such that
$p=\vartheta(x)$). Then by definition of $\mathbf{d}_{h}$, the statement holds.
\end{proof}

\bigskip

By the above result, for any given precision $\eta\in\mathbb{R}^{+}$ we can
approximate the state space $Q_{\tau}=\mathbb{R}^{n}$ of $T_{\tau}(\Sigma)$ by means of the
countable set $Q:=[\mathbb{R}^{n}]_{\eta}$. This choice for $Q$ guarantees that for any $x\in Q_{\tau}$ there exists $q\in Q$ so that $\Vert x-q \Vert \leq \eta/2$. \\
The approximation of the set of labels $L_{\tau}$ of $T_{\tau}(\Sigma)$ is more involved and it requires the notion of reachable set.
We recall that given a forward complete control system $\Sigma=(\mathbb{R}^{n},U\times V,\mathcal{U}%
\times\mathcal{V},f)$, any $\tau\in\mathbb{R}^{+}$ and $x\in\mathbb{R}%
^{n}$, the reachable set of $T_{\tau}(\Sigma)$ with initial condition $x\in Q_{\tau}$ is the set
$\mathcal{R}(\tau,x)$ of endpoints $\mathbf{x}(\tau,x,\mathbf{a,b})$ for any
$\mathbf{a\in}\mathcal{\ }A_{\tau}$ and $\mathbf{b\in}$ $B_{\tau}$, or equivalently:
\begin{equation}
\mathcal{R}(\tau,x):=\left\{y\in Q_{\tau}:
x \rTo_{\tau}^{\mathbf{a,b}}y,\text{ }\mathbf{a}\text{ }\mathbf{\mathbf{\in}}\text{ }A_{\tau},\text{
}\mathbf{b}\in B_{\tau}\right\}.\label{reach1}%
\end{equation}
Moreover, the reachable set of $T_{\tau}(\Sigma)$ with initial condition $x\in Q_{\tau}$ and control
label $\mathbf{a}\in A_{\tau}$ is the set $\mathcal{R}(\tau,x,\mathbf{a})$ of
endpoints $\mathbf{x}(\tau,x,\mathbf{a,b})$ for any $\mathbf{b}\in B_{\tau}$,
i.e.
\begin{equation}
\mathcal{R}(\tau,x,\mathbf{a}):=\left\{y\in Q_{\tau}:x \rTo_{\tau}^{\mathbf{a,b}}y,\text{ }\mathbf{b}\in B_{\tau}\right\}.\label{Reach}%
\end{equation}
The reachable sets in (\ref{reach1}) and (\ref{Reach}) are
well--defined because the control system $\Sigma$ associated with $T_{\tau}(\Sigma)$ is assumed to be forward complete. \\
Given any desired precision $\mu\in \mathbb{R}^{+}$, we approximate $L_{\tau}$ by means of the set \mbox{$L:=A\times B$}, where:
\begin{eqnarray}
A:=%
{\textstyle\bigcup\nolimits_{q\in Q}}
A^{\mu}(q),
&
\hspace{5mm}
&
B:=%
{\textstyle\bigcup\nolimits_{q\in Q}}
{\textstyle\bigcup\nolimits_{\mathbf{a}\in A^{\mu}(q)}}
B^{\mu}(q,\mathbf{a}),
\label{L2}
\end{eqnarray}
and $A^{\mu}(q)$ captures the set of control labels that can be applied at the state $q\in Q$, while 
$B^{\mu}(q,\mathbf{a})$ captures the set of disturbance labels that can be applied at the state $q\in Q$ when the chosen control label is $\mathbf{a}\in A^{\mu}(q)$. The definition of sets $A$ and $B$ in (\ref{L2}) is \textit{asymmetric}. This asymmetry follows from the notion of $A\varepsilon A$ bisimulation relation that we use, where control labels must be chosen 
``robustly'' with respect to the action of disturbance labels (see conditions (ii) and (iii) in Definition \ref{Alt_ASR}). 
%
%
%
%
Given any $\tau\in\mathbb{R}^{+}$, 
 define the following sets:%
\begin{eqnarray}%
&&\texttt{A}_{\mu}(\tau,q):=\{P\in2^{[\mathbb{R}^{n}]_{\mu}}\text{ }|\text{
}\exists\mathbf{a}\in A_{\tau}\text{ s.t. }\mathbf{d}_{h}(P,\mathcal{R}%
(\tau,q,\mathbf{a}))\leq\mu/2\},\label{ccc}
\\
&&\texttt{B}_{\mu}(\tau,q,\mathbf{a}):=\{p\in\lbrack\mathbb{R}^{n}]_{\mu}\text{
}|\text{ }\exists \mathbf{b}\in B_{\tau}\text{ s.t. }\mathbf{d}(p,\mathbf{x}(\tau,q,\mathbf{a},\mathbf{b}))=\Vert
p-\mathbf{x}(\tau,q,\mathbf{a},\mathbf{b})\Vert\leq\mu/2\}.\notag
\end{eqnarray}
Notice that for any $P\in \texttt{A}_{\mu}(\tau,q)$ there exists a (possibly infinite) set of control labels $\mathbf{a}\in A_{\tau}$ so that $\mathbf{d}_{h}(P,\mathcal{R}(\tau,q,\mathbf{a}))\leq\mu/2$. Analogously, for any $p\in \texttt{B}_{\mu}(\tau,q,\mathbf{a})$ there exists a (possibly infinite) set of disturbance labels $\mathbf{b}\in B_{\tau}$ so that $\mathbf{d}(p,\mathbf{x}(\tau,q,\mathbf{a},\mathbf{b}))=\Vert
p-\mathbf{x}(\tau,q,\mathbf{a},\mathbf{b})\Vert\leq\mu/2$. In order to define the sets $A^{\mu}(q)$ and $B^{\mu}(q,\mathbf{a})$ in (\ref{L2}) we consider for any $P\in \texttt{A}_{\mu}(\tau,q)$ only \textit{one} control label $\mathbf{a}\in A_{\tau}$ and respectively for any \mbox{$p\in \texttt{B}_{\mu}(\tau,q,\mathbf{a})$} only \textit{one} disturbance label $\mathbf{b}\in B_{\tau}$, as ``representatives'' of all control labels and all disturbance labels associated with the set $P$ and the vector $p$, respectively. The sets $A^{\mu}(q)$ and $B^{\mu}(q,\mathbf{a})$ will be defined as the collections of these representative control and disturbance labels, respectively. 
The choice of representatives is defined by the functions:
\begin{eqnarray}
\psi_{\mu}^{\tau,q}:\texttt{A}_{\mu}(\tau,q)\rightarrow A_{\tau},
&
\hspace{5mm}
&
\varphi_{\mu}^{\tau,q,\mathbf{a}}:\texttt{B}_{\mu}(\tau,q,\mathbf{a}%
)\rightarrow B_{\tau},
\end{eqnarray}
where:
\begin{itemize}
\item $\psi_{\mu}^{\tau,q}$ associates to any $P\in\texttt{A}_{\mu}(\tau,q)$
one control label\footnote{These control and disturbance labels exist by the definition of the sets $\texttt{A}_{\mu}(\tau,q)$ and $\texttt{B}_{\mu}(\tau,q,\mathbf{a})$.\label{asdf}} $\mathbf{a}=\psi_{\mu}^{\tau,q}(P)$ $\mathbf{\in}$
$A_{\tau}$ so that $\mathbf{d}_{h}(P,\mathcal{R}(\tau,q,\mathbf{a}))\leq\mu/2$;

\item $\varphi_{\mu}^{\tau,q,\mathbf{a}}$ associates to any 
\mbox{$p\in \texttt{B}_{\mu}(\tau,q,\mathbf{a})$} one disturbance label$^{\ref{asdf}}$ \mbox{$\mathbf{b}%
=\varphi_{\mu}^{\tau,q,\mathbf{a}}(p)\mathbf{\in}$ $B_{\tau}$} such that 
$\mathbf{d}(p,\mathbf{x}(\tau,q,\mathbf{a,b}))=\Vert
p-\mathbf{x}(\tau,q,\mathbf{a,b})\Vert\leq\mu/2$.
\end{itemize}
By the above definition, functions $\psi_{\mu}^{\tau,q}$ and $\varphi_{\mu}^{\tau,q,\mathbf{a}}%
$ are not unique. The sets $A^{\mu}(q)$ and $B^{\mu}(q,\mathbf{a})$, appearing in ($\ref{L2}$), can now be defined by:
\begin{eqnarray}
A^{\mu}(q):=\psi_{\mu}^{\tau,q}(\texttt{A}_{\mu}(\tau,q));
&
\hspace{5mm}
&
B^{\mu}(q,\mathbf{a}):=\varphi_{\mu}^{\tau,q,\mathbf{a}}(\texttt{B}_{\mu}(\tau,q,\mathbf{a})).\label{A2B2}
\end{eqnarray}
Since the control system $\Sigma$ is assumed to be forward complete, the reachable sets $\mathcal{R}(\tau,q,\mathbf{a})$, appearing in (\ref{ccc}) are nonempty; hence sets $\texttt{A}_{\mu}(\tau,q)$ and $\texttt{B}_{\mu}(\tau,q,\mathbf{a})$ in (\ref{ccc}) are nonempty and therefore sets $A^{\mu}(q)$ and $B^{\mu}(q,\mathbf{a})$ in (\ref{A2B2}) are nonempty, as well. Furthermore, provided that sets $\texttt{A}_{\mu}(\tau,q)$ and $\texttt{B}_{\mu}(\tau,q,\mathbf{a})$ in (\ref{ccc}) are countable, sets $A^{\mu}(q)$ and $B^{\mu}(q,\mathbf{a})$ in (\ref{A2B2}) are countable, as well. This would guarantee countability of sets of labels $A$ and $B$ in (\ref{L2}) and consequently, countability of the symbolic model in (\ref{T4}). (Note that the set $Q=[\mathbb{R}^{n}]_{\eta}$ in (\ref{T4}) is countable.) In the next step we will state conditions on control systems that guarantee countability of the sets $\texttt{A}_{\mu}(\tau,q)$ and $\texttt{B}_{\mu}(\tau,q,\mathbf{a})$.\\
We now have all the ingredients to define transition system (\ref{T4}).
\begin{definition}
\label{SymbModNLCS}
Given a control system $\Sigma=(\mathbb{R},U\times V,\mathcal{U}\times \mathcal{V},f)$, any $\tau\in\mathbb{R}^{+}$, $\eta
\in\mathbb{R}^{+}$ and $\mu\in\mathbb{R}^{+}$ define the transition
system:
\begin{equation}
T_{\tau,\eta,\mu}(\Sigma):=(Q,L,\rTo,O,H), \label{T3}%
\end{equation}
where:

\begin{itemize}
\item $Q=[\mathbb{R}^{n}]_{\eta};$

\item $L=A\times B$, where:
\begin{eqnarray}
A=
{\textstyle\bigcup\nolimits_{q\in Q}}
A^{\mu}(q),
&
\hspace{5mm}
&
B=%
{\textstyle\bigcup\nolimits_{q\in Q}}
{\textstyle\bigcup\nolimits_{\mathbf{a}\in A^{\mu}(q)}}
B^{\mu}(q,\mathbf{a}),\notag
\end{eqnarray}
and the sets $A^{\mu}(q)$ and $B^{\mu}(q,\mathbf{a})$ are defined in (\ref{A2B2});
\item $q \rTo^{\mathbf{a,b}} p$, if $\mathbf{a}\in
A^{\mu}(q)$, $\mathbf{b}\in B^{\mu}(q,\mathbf{a})$ and $\left\Vert p-\mathbf{x}%
(\tau,q,\mathbf{a,b})\right\Vert \leq\eta/2$;

\item $O=\mathbb{R}^{n}$;

\item $H=\iota:Q\hookrightarrow O$.
\end{itemize}
\end{definition}
Transition system $T_{\tau,\eta,\mu}(\Sigma)$ is metric when we regard
$O=\mathbb{R}^{n}$ as being equipped with the metric $\mathbf{d}(p,q)=\left\Vert
H(p)-H(q)\right\Vert=\left\Vert p-q\right\Vert$; furthermore note that the metric employed for $T_{\tau,\eta,\mu}(\Sigma)$ is the same one used in transition system $T_{\tau}(\Sigma)$.
In the definition of the transition relation $\rTo$ we require $\mathbf{x}(\tau,q,\mathbf{a},\mathbf{b})$ to be in the closed ball $\mathcal{B}_{\eta/2}(p)$. We can instead, require $\mathbf{x}(\tau,q,\mathbf{a},\mathbf{b})$ to be in $\mathcal{B}_{\lambda}(p)$ for any $\lambda\geq\eta/2$. However, we chose $\lambda=\eta/2$ because $\eta/2$ is the smallest value of $\lambda\in\mathbb{R}^{+}$ that ensures $\mathbb{R}^{n}\subseteq{\textstyle\bigcup\nolimits_{q\in\lbrack\mathbb{R}^{n}]_{\eta}}}\mathcal{B}_{\lambda}(q)$. In fact, this choice of $\lambda$ reduces the number of transitions in the definition of the symbolic model in (\ref{T3}).
\bigskip

\textbf{STEP 2.} Transition system $T_{\tau
,\eta,\mu}(\Sigma)$ is not countable, in general, because the set $\texttt{A}%
_{\mu}(\tau,q)$ of (\ref{ccc}) (which is involved in the definition of sets of labels $A$ and $B$) is not so\footnote{Recall that the power set of a countable set is in general not countable (see e.g. \cite{SetTheory}).}. However, if the reachable sets in (\ref{reach1})
associated to $\Sigma$ are bounded, we can guarantee countability of $T_{\tau,\eta,\mu}(\Sigma)$.

\begin{proposition}
\label{Prop1}Consider a forward complete control system $\Sigma=(\mathbb{R}^{n},U\times V,$
$\mathcal{U\times V},f)$ and any $\tau\in\mathbb{R}^{+}$. Suppose that for
any $x\in\mathbb{R}^{n}$ the reachable set\footnote{Note that sets $\mathcal{R}(\tau,x)$ are well--defined because of the forward completeness assumption on the control system.} $\mathcal{R}(\tau,x)$ is bounded.
Then, for any $\eta\in\mathbb{R}^{+}$ and $\mu\in\mathbb{R}^{+}$ the
corresponding transition system $T_{\tau,\eta,\mu}(\Sigma)$ is countable.
\end{proposition}

\begin{proof}
Since for any $\eta\in\mathbb{R}^{+}$ the set of states $Q$ of $T_{\tau,\eta,\mu}(\Sigma)$ is countable, we only
need to show that $A$ and $B$ are countable. 
Given any precision $\mu\in \mathbb{R}^{+}$, for any $q\in Q$ consider the set:%
\[
P(q):=\{p\in\lbrack\mathbb{R}^{n}]_{\mu}\text{
}:\text{ }\exists z\in\mathcal{R}(\tau,q)\text{ s.t. }\Vert p-z\Vert\leq \mu/2\}.
\]
The set $\mathcal{R}(\tau,q)$ is bounded and therefore the
set $P(q)$ is finite. Since for any
$q\in Q$, $\texttt{A}_{\mu}(\tau,q)\subseteq2^{P(q)}$, $\texttt{A}_{\mu}(\tau,q)$ is finite and therefore $A^{\mu}(q)=\psi_{\mu}%
^{\tau,q}(\texttt{A}_{\mu}(\tau,q))$ is finite, as well. Moreover since $A$ is the
union of finite sets $A^{\mu}(q)$ with $q$ ranging in the countable set
$Q=[\mathbb{R}^{n}]_{\eta}$, the set $A$ is countable (see e.g. \cite{SetTheory}). With respect to the set $B$, since for any state $q\in Q$ and any 
$\mathbf{a}\in A^{\mu}(q)$, $\texttt{B}_{\mu}(\tau,q,\mathbf{a})\subseteq \lbrack\mathbb{R}^{n}]_{\mu}$,
the set $\texttt{B}_{\mu}(\tau,q,\mathbf{a})$ is countable and then 
$B^{\mu}(q,\mathbf{a}%
)=\varphi^{\tau,q,\mathbf{a}}_{\mu}(\texttt{B}_{\mu}(\tau,q,\mathbf{a}))$ is countable, as
well. Finally, since $B$ is the union of countable sets $B^{\mu}(q,\mathbf{a})$
with $q$ ranging in the countable set $Q=[\mathbb{R}^{n}]_{\eta}$ and
$\mathbf{a}$ ranging in the finite set $A^{\mu}(q)$, the set $B$ is countable (see e.g. \cite{SetTheory}).
\end{proof}
\bigskip

A direct consequence of the above result is that if the state space of
$\Sigma$ is bounded, which is the case in many realistic situations, the
proposed transition system $T_{\tau,\eta,\mu}(\Sigma)$ is finite.
The following result gives a checkable condition that guarantees countability of $T_{\tau,\eta,\mu
}(\Sigma)$.

\begin{corollary}
\label{Corollary1}Consider a forward complete control system $\Sigma=(\mathbb{R}^{n},U\times V,$
$\mathcal{U\times V},f)$ and suppose that $U \times V$ is compact. 
Then, for any $\tau\in\mathbb{R}^{+}$, $\eta\in\mathbb{R}^{+}$ and $\mu
\in\mathbb{R}^{+}$ the corresponding transition system $T_{\tau,\eta,\mu
}(\Sigma)$ is countable.
\end{corollary}
\begin{proof}
By Proposition 5.1 of \cite{SmoothConverse}, for any $\tau\in\mathbb{R}^{+}$ and $x\in Q_{\tau}$ the
reachable set $\mathcal{R}(\tau,x)$ is bounded. Hence, the result follows by applying Proposition \ref{Prop1}.
\end{proof}
\bigskip

\textbf{STEP 3.}
We can now give the following result, which relates $\delta$--GAS to
existence of (not necessarily countable) symbolic models.

\begin{theorem}
\label{Th_main2}Consider a control system $\Sigma=(\mathbb{R}^{n},U\times
V,\mathcal{U\times V},f)$ and any desired precision $\varepsilon\in
\mathbb{R}^{+}$. If $\Sigma$ is $\delta$--GAS, then for any $\tau\in\mathbb{R}^{+}$, $\mu\in\mathbb{R}^{+}$
and $\eta\in\mathbb{R}^{+}$ satisfying the following condition:%
\begin{equation}
\beta(\varepsilon,\tau)+\mu+\eta/2<\varepsilon, \label{cond2}%
\end{equation}
the corresponding transition system $T_{\tau,\eta,\mu
}(\Sigma)$ is $A\varepsilon A$ bisimilar to $T_{\tau}(\Sigma)$.
\end{theorem}
\bigskip

Before giving the proof of this result we point out that if $\Sigma$ is
$\delta$--GAS, there always exist parameters $\tau\in\mathbb{R}^{+}$,
$\eta\in\mathbb{R}^{+}$ and $\mu\in\mathbb{R}^{+}$ satisfying condition
(\ref{cond2}). In fact, if $\Sigma$ is $\delta$--GAS then there exists a
sufficiently large $\tau\in\mathbb{R}^{+}$ so that $\beta(\varepsilon
,\tau)<\varepsilon$; then by choosing sufficiently small values of $\mu$ and
$\eta$, condition (\ref{cond2}) is fulfilled.

\bigskip

\begin{proof}
Consider the relation $R\subseteq Q_{\tau}\times Q$ defined by $(x,q)\in R$
if and only if $||x-q||\leq\varepsilon$. By construction $R^{-1}(Q)=Q_{\tau}$; by
geometrical considerations on the infinity norm, $Q_{\tau}\subseteq%
{\textstyle\bigcup\nolimits_{p\in\lbrack\mathbb{R}^{n}]_{\eta}}}
\mathcal{B}_{\eta/2}(p)$ and therefore, since by
(\ref{cond2}) $\eta/2<\varepsilon$, we have that $R(Q_{\tau})=Q$. We now
show that $R$ is an $A\varepsilon A$ bisimulation
relation between $T_{\tau}(\Sigma)$ and $T_{\tau,\eta,\mu}(\Sigma)$.

Consider any $(x,q)\in R$. Condition (i) in Definition \ref{Alt_ASR} is
satisfied by the definition of $R$ and of the involved metric transition systems. Let us now show that condition (ii) in
Definition \ref{Alt_ASR} also holds. Since $\delta$--GAS implies forward completeness, reachable sets defined in (\ref{Reach}) are well defined, for any $\tau\in \mathbb{R}^{+}$, $x\in Q_{\tau}$ and $\mathbf{a}\in A_{\tau}$.\\
Consider any $\mathbf{a}_{1}\in A_{\tau}$. Given any $\mu\in \mathbb{R}^{+}$, by Lemma \ref{Lemma1}, there exists
$P\subseteq\lbrack\mathbb{R}^{n}]_{\mu}$ such that:
\begin{equation}
\mathbf{d}_{h}(P,\mathcal{R}(\tau,q,\mathbf{a}_{1}))\leq\mu/2\text{.} \label{a1}%
\end{equation}
By inequality (\ref{a1}), $P\in\texttt{A}_{\mu}(\tau,q)$ and then let $\mathbf{a}_{2}$ be given by\footnote{Note that depending on the choice of function $\psi_{\tau,q}$, which is not unique, $\mathbf{a}_{2}$ can either coincide or not with $\mathbf{a}_{1}$.} $\mathbf{a}_{2}=\psi_{\tau,q}^{\mu}(P)\in
A^{\mu}(q)$. By (\ref{a1}), 
the definition of $\psi_{\tau,q}^{\mu}$ and the properties of $\mathbf{d}_{h}$ we have:%
\begin{equation}
\mathbf{d}_{h}(\mathcal{R}(\tau,q,\mathbf{a}_{1}),\mathcal{R}(\tau,q,\mathbf{a}%
_{2}))\leq\mathbf{d}_{h}(P,\mathcal{R}(\tau,q,\mathbf{a}_{1}))+\mathbf{d}
_{h}(P,\mathcal{R}(\tau,q,\mathbf{a}_{2}))\leq\mu. \label{a2}%
\end{equation}
Consider now any disturbance label\footnote{Existence of such disturbance label is guaranteed by nonemptyness of set $B^{\mu}(q,\mathbf{a}_{2})$.} 
\mbox{$\mathbf{b}_{2}\in B^{\mu}(q,\mathbf{a}_{2})\subset B_{\tau}$} and set
$z=\mathbf{x}(\tau,q,\mathbf{a}_{2},\mathbf{b}_{2})\in \mathcal{R}(\tau,q,\mathbf{a}_{2})$.
By inequality (\ref{a2}) and the definition of $\mathbf{d}_{h}$, 
there exists $z_{1}\in
\overline{\mathcal{R}(\tau,q,\mathbf{a}_{1})}$ such that:%
\begin{equation}
\mathbf{d}(z_{1},z)=\Vert z_{1}-z\Vert\leq\mu. \label{a3}%
\end{equation}
The vector\footnote{The reachable set $\mathcal{R}(\tau,q,\mathbf{a}_{1})$ is in general not closed and therefore inequality (\ref{a2}) does not guarantee the existence of $z_{1}\in
\mathcal{R}(\tau,q,\mathbf{a}_{1})$, satisfying inequality (\ref{a3}). However, by definition of $\mathbf{d}_{h}$, the vector $z_{1}$ is guaranteed to be in the topological closure of the reachable set $\mathcal{R}(\tau,q,\mathbf{a}_{1})$.} $z_{1}$ can be either in $\mathcal{R}(\tau,q,\mathbf{a}_{1})$ or in $\overline{\mathcal{R}(\tau,q,\mathbf{a}_{1})} \,\backslash \, \mathcal{R}(\tau,q,\mathbf{a}_{1})$; in both cases for any $\sigma \in \mathbb{R}^{+}$ there exists $z_{2}\in\mathcal{R}(\tau,q,\mathbf{a}_{1})$ such that:
\begin{equation}
\mathbf{d}(z_{1},z_{2})=\Vert z_{1}-z_{2}\Vert\leq\sigma. \label{a4}%
\end{equation}
(In particular if $z_{1}\in\mathcal{R}(\tau,q,\mathbf{a}_{1})$ one can choose $z_{1}=z_{2}$.)
Choose $\mathbf{b}_{1}\in B_{\tau}$ such that $z_{2}=\mathbf{x}(\tau,q,\mathbf{a}%
_{1},\mathbf{b}_{1})$ (Notice that since $z_{2}\in\mathcal{R}(\tau,q,\mathbf{a}%
_{1})$, such $\mathbf{b}_{1}\in B_{\tau}$ does exist.).
Consider the transition $x\rTo{\mathbf{a}_{1},\mathbf{b}_{1}}_{\tau}y$ in $T_{\tau
}(\Sigma)$. Since $Q_{\tau}\subseteq\bigcup\nolimits_{q'\in [\mathbb{R}^{n}]_{\eta}} \mathcal{B}_{\eta/2}(q')$, there exists $p\in Q=[\mathbb{R}^{n}]_{\eta}$
such that:
\begin{equation}
\mathbf{d}(z,p)=\left\Vert z-p\right\Vert \leq\eta/2. \label{b2}%
\end{equation}
Thus $q\rTo^{\mathbf{a_{2}},\mathbf{b}_{2}}p$ in $T_{\tau,\eta,\mu
}(\Sigma)$. Since
$\Sigma$ is $\delta$--GAS and by (\ref{a4}), (\ref{a3}) and (\ref{b2}) the
following chain of inequalities holds:
\begin{align*}
\Vert y-p\Vert &  =\Vert y-z_{2}+z_{2}-z_{1}+z_{1}-z+z-p\Vert\\
&  \leq\Vert y-z_{2}\Vert+\Vert z_{2}-z_{1}\Vert+\Vert z_{1}-z\Vert+\Vert z-p\Vert\\
&  \leq\beta(||x-q||,\tau)+\Vert z_{2}-z_{1}\Vert+\Vert z_{1}-z\Vert+\Vert z-p\Vert\\
&  \leq\beta(\varepsilon,\tau)+\sigma+\mu+\eta/2.
\end{align*}
By inequality (\ref{cond2}), there exists a sufficiently small value of $\sigma\in\mathbb{R}^{+}$ such that 
$\beta(\varepsilon,\tau)+\sigma+\mu+\eta/2\leq\varepsilon$, and hence $(y,p)\in R$ and condition (ii) in Definition \ref{Alt_ASR} holds. \\
We now show that condition (iii) is also satisfied. Consider any $\mathbf{a}_{2}\in A$; since $A \subset A_{\tau}$, we can choose $\mathbf{a}_{1}%
=\mathbf{a}_{2}\in A_{\tau}$. Consider any $\mathbf{b}_{1}\in B_{\tau}$ and set
$z=\mathbf{x}(\tau,q,\mathbf{a}_{1},\mathbf{b}_{1})$. Since $Q_{\tau}\subseteq
{\textstyle\bigcup\nolimits_{q'\in\lbrack\mathbb{R}^{n}]_{\mu}}}
\mathcal{B}_{\mu/2}(q')$, then there exists
$z_{1}\in\lbrack\mathbb{R}^{n}]_{\mu}$ such that:%
\begin{equation}
\mathbf{d}(z,z_{1})=\Vert z-z_{1}\Vert\leq\mu/2. \label{b3}%
\end{equation}
Furthermore $z\in\mathcal{R}(\tau,q,\mathbf{a}_{1})$ and hence, it is clear that
$z_{1}\in\texttt{B}_{\mu}(\tau,q,\mathbf{a}_{1})$ by definition of $\texttt{B}%
_{\mu}(\tau,q,\mathbf{a}_{1})$. Then let $\mathbf{b}_{2}$ be given by\footnote{Note that depending on the choice of function $\varphi^{\tau,q,\mathbf{a}_{1}}_{\mu}$, which is not unique, $\mathbf{b}_{2}$ can either coincide or not with $\mathbf{b}_{1}$.}
$\mathbf{b}_{2}=\varphi^{\tau,q,\mathbf{a}_{1}}_{\mu}(z_{1})=\varphi^{\tau,q,\mathbf{a}_{2}}_{\mu}(z_{1})\in B^{\mu}%
(q,\mathbf{a}_{2})$. By definition of function $\varphi^{\tau,q,\mathbf{a}_{2}}_{\mu}$ and
by setting $z_{2}=\mathbf{x}(\tau,q,\mathbf{a}_{2},\mathbf{b}_{2})$, it follows that:
\begin{equation}
\mathbf{d}(z_{1},z_{2})=\left\Vert z_{1}-z_{2}\right\Vert \leq\mu/2. \label{b44}%
\end{equation}
Since $Q_{\tau}\subseteq%
{\textstyle\bigcup\nolimits_{q'\in\lbrack\mathbb{R}^{n}]_{\eta}}}
\mathcal{B}_{\eta/2}(q')$, there exists $p\in Q=[\mathbb{R}^{n}]_{\eta}$
such that:
\begin{equation}
\mathbf{d}(z_{2},p)=\left\Vert z_{2}-p\right\Vert \leq\eta/2, \label{b4}%
\end{equation}
and therefore $q\rTo^{\mathbf{a_{2}},\mathbf{b}_{2}}%
p$ in $T_{\tau,\eta,\mu}(\Sigma)$. 
Consider now the transition $x\rTo^{\mathbf{a}_{1},\mathbf{b}_{1}}_{\tau}y$ in $T_{\tau}(\Sigma)$.
Since $\Sigma$ is $\delta$--GAS and by (\ref{b3}), (\ref{b44}), (\ref{b4}) and
(\ref{cond2}), the following chain of inequalities holds:%
\begin{align*}
\Vert y-p\Vert &  =\Vert y-z+z-z_{1}+z_{1}-z_{2}+z_{2}-p\Vert\\
&  \leq\Vert y-z\Vert+\Vert z-z_{1}\Vert+\Vert z_{1}-z_{2}\Vert+\Vert z_{2}-p\Vert\\
&  \leq\beta(\Vert x-q\Vert,\tau)+\Vert z-z_{1}\Vert+\Vert z_{1}-z_{2}\Vert+\Vert z_{2}-p\Vert\\
&  \leq\beta(\varepsilon,\tau)+\mu+\eta/2<\varepsilon.
\end{align*}
Thus $(y,p)\in R$, which completes the proof.
\end{proof}
\bigskip

Finally, by combining Corollary \ref{Corollary1} and Theorem \ref{Th_main2}, the proof of Theorem \ref{Th_main22} holds as a straightforward consequence.

\section{Linear control systems}\label{sec5}

In this section we specialize results of the previous section to the class of linear control systems. The motivation for addressing this special case 
 is twofold: 
\begin{itemize}
\item[(i)] the construction of symbolic models simplifies and can be easily performed; 
\item[(ii)] the proposed symbolic models are $(A_{1},A_{2})$--$(B_{1},B_{2})$--$A\varepsilon A$ bisimilar to linear control systems, while symbolic models defined in (\ref{T3}) are guaranteed to be only $(A_{1},A_{2})$--$A\varepsilon A$ bisimilar to nonlinear control systems.\\
\end{itemize}
A \textit{linear} control system is a control system $\Sigma=(\mathbb{R}^{n},U\times V,\mathcal{U}\times\mathcal{V},f)$, where the vector field $f$ is linear, i.e. for any $x\in \mathbb{R}^{n}$, $u\in U$ and $v\in V$,
\[
f(x,u,v)=\mathbf{A}x+\mathbf{B}u+\mathbf{G}v,
\]
for some matrices $\mathbf{A}$, $\mathbf{B}$ and $\mathbf{G}$ of appropriate dimensions. 
With a slight abuse of notation we say that a linear control system $\Sigma$ is asymptotically stable, when $\Sigma$ with $U \times V = \{0\}\times\{0\}$ is
so. \\
For any given $\tau\in\mathbb{R}^{+}$, consider the following sets:%
\begin{eqnarray}
\mathcal{R}_{A_{\tau}}:=\left\{p\in Q_{\tau}:
0 \rTo_{\tau}^{\mathbf{a,0}}p,\text{ }\mathbf{a}\text{ }\mathbf{\mathbf{\in}}\text{ }A_{\tau}\right\};\notag\\
\mathcal{R}_{B_{\tau}}:=\left\{p\in Q_{\tau}:
0 \rTo_{\tau}^{\mathbf{0,b}}p,\text{ }\mathbf{b}\text{ }\mathbf{\mathbf{\in}}\text{ }B_{\tau}\right\},
\label{RRRRR}
\end{eqnarray}
%
of reachable states of $T_{\tau}(\Sigma)$ from the origin $0$ by means of any control
label $\mathbf{a}\in A_{\tau}$ and identically null disturbance label $\mathbf{0}$ and,
respectively, by means of any disturbance label $\mathbf{b}\in B_{\tau}$
 and identically null control label $\mathbf{0}$. Notice that sets in (\ref{RRRRR}) are well--defined since label curves 
$\mathbf{a}$ and $\mathbf{b}$ are locally essentially bounded. 
We can now propose the following symbolic models for linear systems.
\begin{definition}
Given a linear control system $\Sigma=(\mathbb{R}^{n},U\times V,\mathcal{U\times V},f)$ and any $\tau\in\mathbb{R}^{+}$, $\eta
\in\mathbb{R}^{+}$ and $\mu\in\mathbb{R}^{+}$,
define the following transition system:
\begin{equation}
\label{T_lin}
T_{\tau,\eta,\mu}(\Sigma):=(Q,A\times B,\rTo,O,H),
\end{equation}
where:

\begin{itemize}
\item $Q=[\mathbb{R}^{n}]_{\eta};$

\item $A$ is a subset of $[\mathbb{R}^{n}]_{\mu}$ for which $\mathbf{d}
_{h}(A,\mathcal{R}_{A_{\tau}})\leq\mu/2$;

\item $B$ is a subset of $[\mathbb{R}^{n}]_{\mu}$ for which $\mathbf{d}
_{h}(B,\mathcal{R}_{B_{\tau}})\leq\mu/2$;

\item $q \rTo^{a,b} p$, if the following inequality is satisfied 
\begin{equation}
\left\Vert \mathbf{x}(\tau,q,0,0)+a+b-p\right\Vert \leq\eta/2;
\label{ffff}
\end{equation}

\item $O=\mathbb{R}^{n}$;

\item $H=\iota:Q\hookrightarrow O$.
\end{itemize}
\end{definition}
By Lemma \ref{Lemma1}, sets of labels $A$ and
$B$ do exist; moreover they are \textit{countable}. Hence, transition system $T_{\tau,\eta,\mu}(\Sigma)$ is countable, as well. Notice that in this case we do not need to require the input set $U\times V$ to be compact for ensuring countability of $T_{\tau,\eta,\mu}(\Sigma)$, whereas in the case of nonlinear control systems it was indeed required (see Corollary \ref{Corollary1}). Furthermore, transition system $T_{\tau,\eta,\mu}(\Sigma)$ of (\ref{T_lin}) can be easily constructed.
The construction of $T_{\tau,\eta,\mu}(\Sigma)$ relies on the computation of the reachable sets in (\ref{RRRRR}). The exact computation of those sets
is in general hard. 
However,
there are
several results available in the literature, that propose approximations of 
reachable sets for linear control systems (e.g. \cite{Varaiya:98,Girard_HSCC05,Elipsoid,Tomlin_Reach} and the references therein). 
For example following \cite{Varaiya:98,Girard_HSCC05}, if $U\times V$ is compact and $\mathcal{U}\times \mathcal{V}$ is the class of all measurable and essentially bounded functions taking values in $U\times V$, given any precisions $e_{A_{\tau}},e_{B_{\tau}}\in\mathbb{R}^{+}$, 
it is possible to compute a pair of polytopes $P_{e_{A_{\tau}}}(\mathcal{R}_{A_{\tau}})$, $P_{e_{B_{\tau}}}(\mathcal{R}_{B_{\tau}})$ so that: 
\[%
\begin{array}
[c]{ccc}%
\mathbf{d}_{h}(P_{e_{A_{\tau}}}(\mathcal{R}_{A_{\tau}}),\mathcal{R}_{A_{\tau}})\leq e_{A_{\tau}}, &
& \mathbf{d}_{h}(P_{e_{B_{\tau}}}(\mathcal{R}_{B_{\tau}}),\mathcal{R}_{B_{\tau}}%
)\leq e_{B_{\tau}}.
\end{array}
\]
Once sets $P_{e_{A_{\tau}}}%
(\mathcal{R}_{A_{\tau}})$ and $P_{e_{B_{\tau}}}(\mathcal{R}_{B_{\tau}})$ are known\footnote{
The interested reader can refer to \cite{Varaiya:98} for an analysis of the
computational effort, required for computing polytopic approximations of
reachable sets for linear control systems.}, the computation of sets of labels $A$ and $B$ can be performed, as well. In fact, sets $A$ and $B$ can be computed on the basis of the approximating sets $P_{e_{A_{\tau}}}(\mathcal{R}_{A_{\tau}})$ and $P_{e_{B_{\tau}}}(\mathcal{R}_{B_{\tau}})$, rather than on the basis of reachable sets $\mathcal{R}_{A_{\tau}}$ and $\mathcal{R}_{B_{\tau}}$. The numerical errors $e_{A_{\tau}}$ and $e_{B_{\tau}}$ can be incorporated in the symbolic model in (\ref{T_lin}), by replacing inequality (\ref{ffff}) by 
\begin{equation}
\left\Vert \mathbf{x}(\tau,q,0,0)+a+b-p\right\Vert \leq\eta/2+e_{A_{\tau}}+e_{B_{\tau}}.
\label{zzz}
\end{equation}
Moreover, 
%
since the norm in inequality (\ref{ffff})
is the infinity norm, the construction of transition relation $\rTo$ can be performed,
by using standard techniques available in the literature on linear
matrix inequalities \cite{LMI}. Indeed, by defining 
\mbox{$\mathcal{B}_{1}(0):=\{x\in \mathbb{R}^{n}:\texttt{M}x\leq\texttt{m}\}$}
 for some matrix $\texttt{M}$ and vector $\texttt{m}$, condition (\ref{ffff}) becomes:
\[
\texttt{M}(e^{A\tau}q+a+b-p)\leq (\eta/2) \,\texttt{m}.
\]

Note that transition system $T_{\tau,\eta,\mu}(\Sigma)$ in (\ref{T_lin}) differs from the one proposed in (\ref{T3}) for nonlinear control systems, (only) in the definition of the sets of labels $A$ and $B$. (It is easy to see that the definition of the transition relation in the two symbolic models is substantially equivalent.) In particular, the set of disturbance labels $B$ is defined \textit{independently} from the set of control labels $A$. This feature is in fact a direct consequence of the linearity assumption, resulting in the so-called 
superposition principle. We can now give the following result.

\begin{theorem}
\label{Th_lin}
Consider a linear control system $\Sigma=(\mathbb{R}^{n},U\times
V,\mathcal{U\times V},f)$ and any
desired precision $\varepsilon\in\mathbb{R}^{+}$. If $\Sigma$ is asymptotically stable then for any $\tau\in\mathbb{R}^{+}$,
$\mu\in\mathbb{R}^{+}$ and $\eta\in\mathbb{R}^{+}$ satisfying the following
condition:%
\begin{equation}
\left\Vert e^{\mathbf{A}\tau}\right\Vert \varepsilon+\mu+\eta/2
<\varepsilon, %
\label{condlin}
\end{equation}
the corresponding transition system $T_{\tau,\eta,\mu}(\Sigma)$ is $A\varepsilon A$ bisimilar to $T_{\tau}(\Sigma)$.
\end{theorem}

\begin{proof}
Consider the relation $R\subseteq Q_{\tau}\times Q$ defined by $(x,q)\in R$
if and only if $||x-q||\leq\varepsilon$. As shown in the proof of Theorem \ref{Th_main2}, $R^{-1}(Q)=Q_{\tau}$, $R(Q_{\tau})=Q$ and $R$ satisfies condition (i) in Definition \ref{Alt_ASR}. We now show that $R$ satisfies also conditions (ii) and (iii).

Consider any $(x,q)\in R$, 
any $\mathbf{a}_{1}\in A_{\tau}$ and choose $a_{2}\in A$ such that:
\begin{equation}
\left\Vert a_{2}-\mathbf{x}(\tau,0,\mathbf{a}_{1},\mathbf{0})\right\Vert
\leq\mu/2. \label{c4}%
\end{equation}
(Such control label $a_{2}$ exist because the set $A$ is closed.)
Consider any $b_{2}\in B$. By the definition of $B$ and of $\mathbf{d}_{h}$, 
there exists $b_{3}\in
\overline{\mathcal{R}_{B_{\tau}}}$ such that:%
\begin{equation}
\mathbf{d}(b_{2},b_{3})=\Vert b_{2}-b_{3}\Vert\leq\mu/2. \label{bb3}%
\end{equation}
The vector $b_{3}$ can be either in $\mathcal{R}_{B_{\tau}}$ or in $\overline{\mathcal{R}_{B_{\tau}}} \,\backslash \, \mathcal{R}_{B_{\tau}}$; in both cases for any $\sigma \in \mathbb{R}^{+}$ there exists $b_{4}\in \mathcal{R}_{B_{\tau}}$ such that:
\begin{equation}
\left\Vert b_{3}-b_{4}
\right\Vert\leq\sigma.
\label{bbb2}
\end{equation}
Choose $\mathbf{b}_{1}\in B_{\tau}$ such that $b_{4}=\mathbf{x}(\tau,0,\mathbf{0},\mathbf{b}_{1})$ and
consider the transition $x\rTo^{\mathbf{a}_{1}\mathbf{,b}_{1}%
}_{\tau}y$ in $T_{\tau}(\Sigma)$. 
Set $z=\mathbf{x}(\tau,q,\mathbf{0},\mathbf{0})+a_{2}+b_{2}\in Q_{\tau}$; since $Q_{\tau}\subseteq
{\textstyle\bigcup\nolimits_{q'\in\lbrack\mathbb{R}^{n}]_{\eta}}}
\mathcal{B}_{\eta/2}(q')$, there exists $p\in Q=[\mathbb{R}^{n}]_{\eta}$
such that:
\begin{equation}
\left\Vert z-p\right\Vert \leq\eta/2. \label{c2}%
\end{equation}
Thus $q\rTo^{a_{2},b_{2}%
}p$ in $T_{\tau,\eta,\mu}(\Sigma)$. By inequalities (\ref{c2}), (\ref{c4}), 
(\ref{bbb2}) and (\ref{condlin}), the following chain of inequalities holds:%
\begin{eqnarray}
\Vert y-p\Vert &&  =\Vert y-z+z-p\Vert\notag\\
&&  \leq\Vert y-z\Vert+\Vert z-p\notag\Vert\\
&&  \leq\Vert\mathbf{x}(\tau,x,\mathbf{a}_{1},\mathbf{b}_{1})-(\mathbf{x}%
(\tau,q,\mathbf{0},\mathbf{0})+a_{2}+b_{2})\Vert+\eta/2\notag\\
&&  =\Vert\mathbf{x}(\tau,x-q,\mathbf{0},\mathbf{0})+\mathbf{x}(\tau,0,\mathbf{a}_{1}%
,\mathbf{0})-a_{2}+\mathbf{x}(\tau,0,\mathbf{0},\mathbf{b}_{1})-b_{2}+b_{3}-b_{3}\Vert+\eta/2\notag\\
&&  \leq\Vert e^{\mathbf{A}\tau}(x-q)\Vert+\Vert\mathbf{x}(\tau,0,\mathbf{a}%
_{1},\mathbf{0})-a_{2}\Vert+\Vert b_{4}-b_{3}\Vert+\Vert b_{3}-b_{2} \Vert +
\eta/2\notag\\
&&  \leq\Vert e^{\mathbf{A}\tau}\Vert\varepsilon+\mu/2+\sigma+\mu/2+\eta/2\text{.}%
\label{aaa}
\end{eqnarray}
By inequality (\ref{condlin}), there exists a sufficiently small value of $\sigma\in\mathbb{R}^{+}$ such that 
\mbox{$\beta(\varepsilon,\tau)+\sigma+\mu+\eta/2\leq\varepsilon$}, and hence $(y,p)\in R$ and condition (ii) in Definition \ref{Alt_ASR} holds. \\
Condition (iii) of Definition \ref{Alt_ASR} can be shown by using same arguments and therefore it is omitted.
\end{proof}
\bigskip

We stress that conditions of Theorem \ref{Th_lin} are conceptually equivalent to conditions of Theorem \ref{Th_main2}. In fact $\delta$--GAS for linear control systems is equivalent to the asymptotic stability of matrix $\mathbf{A}$.
Furthermore, it is well known (e.g. \cite{ISS}) that for linear control systems, function $\beta$ appearing in inequality (\ref{deltaUGAS}), can be chosen as $\beta(\left\Vert x_{1}-x_{2}\right\Vert ,t)=\left\Vert e^{\mathbf{A}t}\right\Vert \left\Vert x_{1}-x_{2}\right\Vert$, and therefore condition (\ref{cond2}) boils down in this case to condition (\ref{condlin}). 
Although assumptions on Theorem \ref{Th_lin} and Theorem \ref{Th_main2} coincide for the class of linear control systems, we stress that Theorem \ref{Th_lin} and Theorem \ref{Th_main2} relate $T_{\tau}(\Sigma)$ to symbolic models in (\ref{T3}) and (\ref{T_lin}), respectively. While the construction of symbolic model in (\ref{T3}) is hard in general, as pointed out before symbolic model in (\ref{T_lin}) can be easily constructed. Furthermore, Theorem \ref{Th_lin} can be extended to a more general result that we state hereafter.
\begin{corollary}
Consider a linear control system $\Sigma=(\mathbb{R}^{n},U\times
V,\mathcal{U\times V},f)$ and any
desired precision $\varepsilon\in\mathbb{R}^{+}$. If $\Sigma$ is asymptotically stable then for any \mbox{$\tau\in\mathbb{R}^{+}$},
\mbox{$\mu\in\mathbb{R}^{+}$} and $\eta\in\mathbb{R}^{+}$ satisfying condition (\ref{condlin}), transition system  $T_{\tau,\eta,\mu}(\Sigma)$ is 
\mbox{$(A_{\tau},A)$--} $(B_{\tau},B)$--$A\varepsilon A$ bisimilar to 
$T_{\tau}(\Sigma)$.
\label{cor}
\end{corollary}
\bigskip

The proof of the result above is a straightforward consequence of the symmetric definition of set of labels $A\times B$ in transition system $T_{\tau,\eta,\mu}(\Sigma)$ and is therefore omitted. Note that the same reasoning does not apply to the general case of nonlinear control systems. The result above is important from a game theory point of view. Suppose that the goal is to find a symbolic model for an infinite state game, the arena of which, is given by a linear control system. Then, Corollary \ref{cor} provides a way of finding a symbolic model that can be used at the same time to design strategies both for the protagonist and for the antagonist of the game.

\section{Illustrative Example}\label{sec6}
In this section we illustrate the results of the previous section in the context of direct current motors. Consider the simplified model of a direct current motor:
\begin{eqnarray}
\Sigma=\left\{
\begin{array}
[c]{l}%
\dot{x}=\mathbf{A}x+\mathbf{B}u+\mathbf{G}v,\\
x\in X, u\in U, v\in V,
\end{array}
\right.
\label{sysexample}
\end{eqnarray}
where $x=(x_{1}$ $x_{2})'$, $x_{1}$ is the current, $x_{2}$ is the angular velocity, $u$ is the applied voltage, $v$ is the load torque disturbance and: 
\begin{align}
\mathbf{A}  &  =\left[
\begin{array}
[c]{rr}%
-R/L & -k_{b}k_{m}/L\\
1/J & -k_{f}/J
\end{array}
\right]  ;\mathbf{B}=\left[
\begin{array}
[c]{c}%
k_{m}/L\\
0
\end{array}
\right]  ;\mathbf{G}=\left[
\begin{array}
[c]{c}%
0\notag\\
1/J
\end{array}
\right],  
\end{align}
where $R=2$ and $L=0.5$ are the resistance and the inductance associated with the armature of the direct current motor; $k_{b}=0.1$ is the back electromagnetic force; $k_{m}=0.1$ is the torque constant; $k_{f}=0.2$ is the viscous friction constant and $J=0.4$ is the inertia momentum. We suppose that:
\[
X=[0,0.6]\times [0,0.6];\text{     }U=[0.3,0.7];\text{     }V=[-0.02,0.02].\notag
\]
All variables and constants appearing in system (\ref{sysexample}) are expressed in the International System. The control problem that we focus on is the one of \textit{disturbance rejection} and it consists in finding a (memoryless) control strategy $\mathbf{u}$, so that for any initial condition $x\in X$ and any disturbance $\mathbf{v}\in \mathcal{V}$, the corresponding angular velocity $x_{2}$ at time $t=5$ is above $0.1$, or equivalently:
\begin{equation}
\mathbf{x}(5,x,\mathbf{u},\mathbf{v})\in X^{*}:=[0,0.6]\times [0.1,0.6].
\label{DistPb}
\end{equation}
Since the system in (\ref{sysexample}) is asymptotically stable, we can apply Theorem \ref{Th_lin}. Set the precision $\varepsilon=0.5$ and $\tau=5$. 
By choosing $\mu=0.3$ and $\eta=0.15$, inequality (\ref{condlin}) is satisfied and therefore the transition system $T_{5,0.3,0.15}(\Sigma)$ defined in (\ref{T_lin}), is $A\varepsilon A$ bisimilar to $T_{5}(\Sigma)$ with $\varepsilon=0.5$.\\
\begin{figure}[ptb]
\begin{center}
\includegraphics[angle=90, scale=0.55]{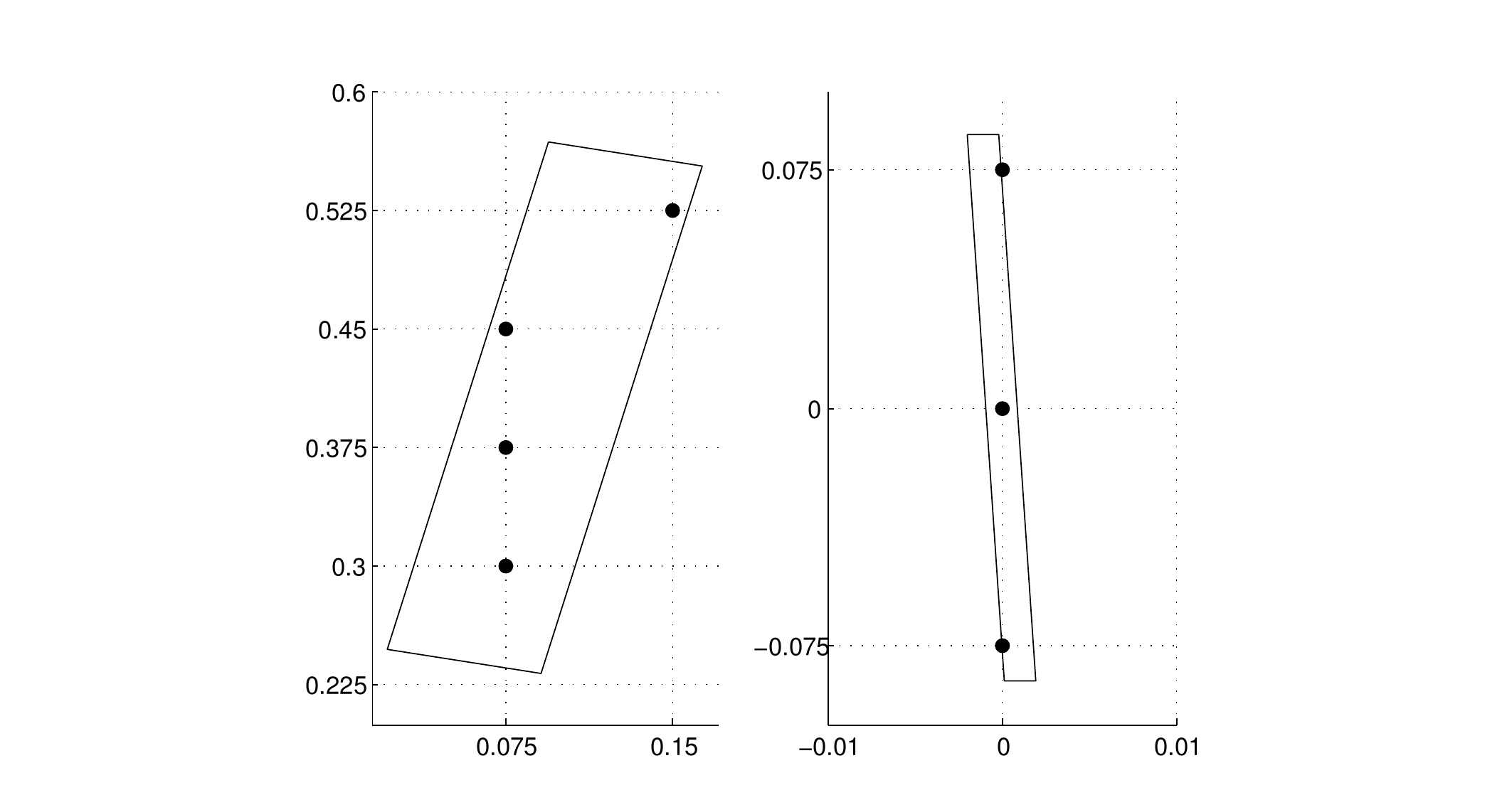}
\end{center}
\caption{Left panel: Outer approximation $P_{e_{A_{\tau}}}(\mathcal{R}_{A_{\tau}})$ of reachable set $\mathcal{R}_{A_{\tau}}$ and control labels set $A$ (black spots). Right panel: Outer approximation $P_{e_{B_{\tau}}}(\mathcal{R}_{B_{\tau}})$ of reachable set $\mathcal{R}_{B_{\tau}}$ and disturbance labels set $B$ (black spots).}%
\label{reachUV}%
\end{figure}
The construction of transition system $T_{5,0.3,0.15}(\Sigma)$ requires the computation of the reachable sets $\mathcal{R}_{A_{\tau}}$ and $\mathcal{R}_{B_{\tau}}$, as defined in (\ref{RRRRR}). By using results in \cite{Girard_HSCC05} and the toolbox MATISSE \cite{MATISSE}, it is possible to compute the following polytopic outer approximations $P_{e_{A_{\tau}}}(\mathcal{R}_{A_{\tau}})$ and $P_{e_{B_{\tau}}}(\mathcal{R}_{B_{\tau}})$ of $\mathcal{R}_{A_{\tau}}$ and $\mathcal{R}_{B_{\tau}}$, respectively:
\begin{eqnarray}
&&
P_{e_{A_{\tau}}}(\mathcal{R}_{A_{\tau}})=conv \left(
\left(
\begin{array}
[c]{c}%
0.1627\\
0.5524
\end{array}
\right),
\left(
\begin{array}
[c]{c}%
0.0908\\
0.2320
\end{array}
\right),
\left(
\begin{array}
[c]{c}%
0.0220\\
0.2474
\end{array}
\right),
\left(
\begin{array}
[c]{c}%
0.0939\\
0.5678
\end{array}
\right)\right);
\label{R_example}
\\
&&
P_{e_{B_{\tau}}}(\mathcal{R}_{B_{\tau}})=conv \left(
\left(
\begin{array}
[c]{c}%
-0.0002\\
0.0862
\end{array}
\right),
\left(
\begin{array}
[c]{c}%
0.0002\\
-0.0862
\end{array}
\right),
\left(
\begin{array}
[c]{c}%
0.0002\\
-0.0862
\end{array}
\right),
\left(
\begin{array}
[c]{c}%
-0.0002\\
0.0862
\end{array}
\right)\right).
\notag
\end{eqnarray}
Numerical errors $e_{A_{\tau}}$ and $e_{B_{\tau}}$ for the sets in (\ref{R_example}), can be evaluated by using Lemma 1 of \cite{Girard_HSCC05}, resulting in $e_{A_{\tau}}=3.0453\,\cdot \,10^{-6}$ and $e_{B_{\tau}}=3.8067\,\cdot\,10^{-5}$. Since $\varepsilon>>\max\{e_{A_{\tau}},e_{B_{\tau}}\}$ we will neglect errors $e_{A_{\tau}}$ and $e_{B_{\tau}}$ in the following developments. (However, as pointed out in the previous section, numerical errors $e_{A_{\tau}}$ and $e_{B_{\tau}}$ could be incorporated in the symbolic model (\ref{T_lin}), by replacing inequality (\ref{ffff}) by (\ref{zzz}).)
On the basis of the sets in (\ref{R_example}) we can compute the sets of labels $A$ and $B$ of transition system (\ref{T_lin}), as shown in Figure \ref{reachUV}.
\begin{table}
\footnotesize
\begin{center}
\begin{tabular}{|c|c|c|c|c|c|c|c|c|c|}
\hline
$q\rTo^{a,b}p$&$q_{1}$ & $q_{2}$ & $q_{3}$ & $q_{4}$ & $q_{5}$ & $q_{6}$ & $q_{7}$ & $q_{8}$ & $q_{9}$ \\
\hline
$a_{1},b_{1}$ & $q_{6}$ & $q_{3}$ & $q_{3}$ & $q_{3}$ & $q_{3}$ & $q_{3}$ & -- & -- & --\\
$a_{1},b_{2}$ & $q_{6}$ & $q_{3}$ & $q_{3}$ & $q_{3}$ & $q_{3}$ & $q_{3}$ & $q_{3}$ & $q_{3}$ & $q_{3}$\\
$a_{1},b_{3}$ & $q_{6}$ & $q_{3}$ & $q_{3}$ & $q_{3}$ & $q_{3}$ & $q_{3}$ & $q_{3}$ & $q_{3}$ & $q_{3}$\\
\hline
$a_{2},b_{1}$ & $q_{3}$ & $q_{3}$ & $q_{3}$ & $q_{3}$ & $q_{3}$ & $q_{3}$ & $q_{3}$ & $q_{3}$ & $q_{3}$ \\
$a_{2},b_{2}$ & $q_{3}$ & $q_{3}$ & $q_{3}$ & $q_{3}$ & $q_{3}$ & $q_{3}$ & $q_{3}$ & $q_{3}$ & $q_{3}$\\
$a_{2},b_{3}$ & $q_{2}$ & $q_{2}$ & $q_{2}$ & $q_{3}$ & $q_{3}$ & $q_{3}$ & $q_{3}$ & $q_{3}$ & $q_{3}$\\
\hline
$a_{3},b_{1}$ & $q_{2}$ & $q_{3}$ & $q_{3}$ & $q_{3}$ & $q_{3}$ & $q_{3}$ & $q_{3}$ & $q_{3}$ & $q_{3}$\\
$a_{3},b_{2}$ & $q_{2}$ & $q_{2}$ & $q_{2}$ & $q_{3}$ & $q_{3}$ & $q_{3}$ & $q_{3}$ & $q_{3}$ & $q_{3}$\\
$a_{3},b_{3}$ & $q_{2}$ & $q_{2}$ & $q_{2}$ & $q_{2}$ & $q_{2}$ & $q_{2}$ & $q_{3}$ & $q_{3}$ & $q_{3}$\\
\hline
$a_{4},b_{1}$ & $q_{2}$ & $q_{2}$ & $q_{2}$ & $q_{3}$ & $q_{3}$ & $q_{3}$ & $q_{3}$ & $q_{3}$ & $q_{3}$\\
$a_{4},b_{2}$ & $q_{2}$ & $q_{2}$ & $q_{2}$ & $q_{2}$ & $q_{2}$ & $q_{2}$ & $q_{3}$ & $q_{3}$ & $q_{3}$\\
$a_{4},b_{3}$ & $q_{2}$ & $q_{2}$ & $q_{2}$ & $q_{2}$ & $q_{2}$ & $q_{2}$ & $q_{2}$ & $q_{2}$ & $q_{2}$\\
\hline
\end{tabular}
\end{center}
\caption{Transition relation of transition system $T_{5,0.3,0.15}(\Sigma)$ defined in (\ref{Texample}). Entry $q_{6}$ corresponding to the second row and second column means that there is a transition from $q_{1}$ to $q_{6}$ labeled by $a_{1},b_{1}$. Entries ``--'' correspond to transitions that end up outside the set $X$.}
\label{TRExample}
\end{table}
The resulting symbolic model:
\begin{equation}
\label{Texample}
T_{5,0.3,0.15}(\Sigma):=(Q,A\times B,\rTo,O,H),
\end{equation}
of (\ref{T_lin}) is given by:

\begin{itemize}
\item $Q=\{q_{1},q_{2},q_{3},q_{4},q_{5},q_{6},q_{7},q_{8},q_{9}\}$, where $q_{1}=(0,0)'$, $q_{2}=(0,\eta)'$, \mbox{$q_{3}=(0,2\eta)'$,} $q_{4}=(\eta,0)'$, $q_{5}=(\eta,\eta)'$, $q_{6}=(\eta,2\eta)'$, $q_{7}=(2\eta,0)'$, $q_{8}=(2\eta,\eta)'$, \mbox{$q_{9}=(2\eta,2\eta)'$;}
\item $A=\{a_{1},a_{2},a_{3},a_{4}\}$, where $a_{1}=(0.1500,0.5250)'$, $a_{2}=(0.0750,0.4500)'$, $a_{3}=(0.0750,0.3750)'$ and $a_{4}=(0.0750,0.3000)'$;
\item $B=\{b_{1},b_{2},b_{3}\}$, where $b_{1}=(0,0.0750)'$, $b_{2}=(0,0)'$ and $b_{3}=(0,-0.0750)'$;

\item $q \rTo^{a,b} p$ is defined in Table \ref{TRExample};

\item $O=\mathbb{R}^{2}$;

\item $H=\iota:Q\hookrightarrow O$.
\end{itemize}
The above symbolic model is depicted in Figure \ref{fig1}. By Theorem \ref{Th_lin} transition systems $T_{5,0.3,0.15}(\Sigma)$ and $T_{5}(\Sigma)$ are $A\varepsilon A$ bisimilar with $\varepsilon=0.5$; furthermore it is easy to see that:
\[
X^{*}=\mathcal{B}_{0.5}(q_3)\cup\mathcal{B}_{0.5}(q_6)\cup\mathcal{B}_{0.5}(q_9).
\]
Hence, the disturbance rejection problem can be solved on the symbolic model in (\ref{Texample}), by finding for any state $q\in Q$, the set $U^{*}(q)$ of all control labels $a\in A$ so that $q\rTo^{a,b}p\in \{q_{3},q_{6},q_{9}\}$ for any disturbance label $b\in B$. A simple inspection of Table \ref{TRExample} provides the following solution to the control problem on the symbolic model:
\begin{eqnarray}
&&
U^{*}(q_{1})=U^{*}(q_{2})=U^{*}(q_{3})=\{a_{1}\}; \notag \\
&&
U^{*}(q_{4})=U^{*}(q_{5})=U^{*}(q_{6})=\{a_{1},a_{2}\}; \notag \\
&&
U^{*}(q_{7})=U^{*}(q_{8})=U^{*}(q_{9})=\{a_{2},a_{3}\}.  
\label{Uopt}
\end{eqnarray}
On the basis of control labels in (\ref{Uopt}) it is possible to synthesize controllers for solving the disturbance rejection problem on the original linear control system (\ref{sysexample}). Indeed, since the system in (\ref{sysexample}) is controllable, by using standard results in linear control theory (see e.g. \cite{Sontag}) for any control label $a\in\{a_{1},a_{2},a_{3}\}\subset\mathcal{R}^{\tau}_{U}$ in (\ref{Uopt}) it is possible to compute a control input $\mathbf{u}\in\mathcal{U}$ so that:
\[
a=\int_{0}^{\tau}e^{A(\tau-t)}\mathbf{B}\mathbf{u}(t)dt.
\]
By definition of the transition system in (\ref{Texample}), the obtained control inputs solve the disturbance rejection problem on the original system in (\ref{sysexample}).

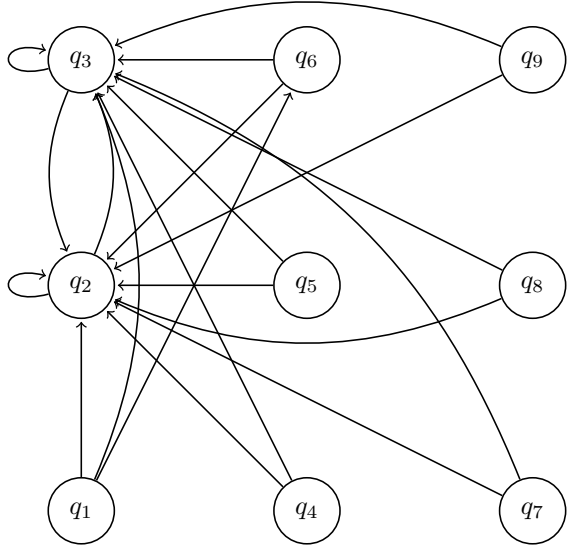
\begin{figure}[ptb]
\begin{center}
\begin{tikzpicture}[->,shorten >=1pt,%
auto,node distance=3cm,semithick,
inner sep=1.5pt,bend angle=23]
\node[state] (q_1) {$q_{3}$};
\node[state] (q_2) [ right of=q_1] {$q_{6}$};
\node[state] (q_3) [ right of=q_2] {$q_{9}$};
\node[state] (q_4) [ below of=q_1] {$q_{2}$};
\node[state] (q_5) [ right of=q_4] {$q_{5}$};
\node[state] (q_6) [ right of=q_5] {$q_{8}$};
\node[state] (q_7) [ below of=q_4] {$q_{1}$};
\node[state] (q_8) [ right of=q_7] {$q_{4}$};
\node[state] (q_9) [ right of=q_8] {$q_{7}$};
\path[->] (q_1) edge [bend right] node {} (q_4)
edge [loop left] node {} ()
(q_2) edge node {} (q_1)
edge node {} (q_4)
(q_3) edge node {} (q_4)
edge [bend right] node {} (q_1)
(q_4) edge [bend right] node {} (q_1)
edge [loop left] node {} ()
(q_5) edge node {} (q_1)
edge node {} (q_4)
(q_6) edge [bend left] node {} (q_4)
edge node {} (q_1)
(q_7) edge [bend right] node {} (q_1)
edge node {} (q_2)
edge node {} (q_4)
(q_8) edge node {} (q_1)
edge node {} (q_4)
(q_9) edge [bend right] node {} (q_1)
edge node {} (q_4);
\end{tikzpicture}
\end{center}
\caption{Symbolic model $T_{5,0.3,0.15}(\Sigma)$ of (\ref{Texample}) associated with the linear control system $\Sigma$ defined in (\ref{sysexample}). An arrow from a state $q$ to a state $p$ means that there exists \textit{at least} a pair $(a,b)\in A\times B$ so that $\mathbf{x}(5,q,0,0)+a+b$ is in the closed ball $\mathcal{B}_{0.30/2}(p)$.}
\label{fig1}%
\end{figure}

\section{Discussion\label{sec7}}
In this paper we showed existence of symbolic models that are $A\varepsilon A$ bisimilar to $\delta$--GAS nonlinear control systems with disturbances. Moreover, the parameter $\varepsilon$ describing the precision, can be chosen as small as desired. For the special class of (asymptotically stable) linear control systems the resulting symbolic models are not only easily computable but, they are also $(A_{1},A_{2})$--$(B_{1},B_{2})$--$A\varepsilon A$ bisimilar to the original systems.\\
The results of this paper generalize the work in \cite{pola-2007} to control systems in presence of disturbances (compare Theorem \ref{Th_main2} and Theorem 4.2 of \cite{pola-2007}). While Theorem 4.2 of \cite{pola-2007} states existence of symbolic models that are approximately bisimilar (in the sense of Definition \ref{ASR_bis}) to $\delta$--GAS control systems, Theorem \ref{Th_main2} shows existence of symbolic models that are $A\varepsilon A$ bisimilar to control systems influenced by disturbances.
%
%
%
As pointed out in Section \ref{Sec32}, the results of \cite{pola-2007} cannot directly be applied to the case of control systems with disturbances. Indeed as Example \ref{exa} shows, the symbolic model (4.4) of \cite{pola-2007} do not capture the different role played by the control inputs and by the disturbance inputs. As a consequence, control strategies synthesized on the symbolic model of \cite{pola-2007} \textit{cannot be transfered} to the original system.
The same observation applies to the results in \cite{BisimSchaft}. As the focus of \cite{BisimSchaft} was the reduction of control systems and not control design, the employed notion of bisimulation was a variation the one of Milner \cite{Milner} and Park \cite{Park}. However, as in the case of the results in \cite{pola-2007}, the notion of bisimulation in \cite{BisimSchaft} cannot be used for control design.
\\
This paper also shares similar ideas with \cite{PolaTabuadaCDC07b}. The work in \cite{PolaTabuadaCDC07b} proposes symbolic models for linear control systems with disturbances. The approximation notion employed in \cite{PolaTabuadaCDC07b} is $A\varepsilon A$ simulation (one--sided version of $A\varepsilon A$ bisimulation). The results in this paper extend the ones in \cite{PolaTabuadaCDC07b} by: 
\begin{itemize}
\item[(i)] enlarging the class of control systems from linear to nonlinear;
\item[(ii)] enlarging the class of control inputs from piecewise constant to measurable and locally essentially bounded;
\item[(iii)] generalizing results from simulation to bisimulation.
\end{itemize}
In particular by (iii), the symbolic model in Definition \ref{ffff} provides a more accurate description of the control system than the one proposed in \cite{PolaTabuadaCDC07b}. This is essential for controller synthesis since, if a controller fails to exist for the symbolic model in \cite{PolaTabuadaCDC07b}, nothing can be concluded regarding the existence of a controller for the original control system. Our results guarantee, instead, that given a control system and a specification, a controller exists for the original model if and only if a controller exists for the symbolic model, up to the resolution $\varepsilon$.\\
Future work will concentrate on constructive techniques to obtain the symbolic models whose existence was shown in this paper.
\bigskip

\textsc{Acknowledgment.} The authors would like to thank Antoine Girard (Universit\'e Joseph Fourier, France) for stimulating discussions on the topic of this paper.

\bigskip

\bibliographystyle{alpha}
\bibliography{biblio}

\end{document}